\documentclass[smallextended]{svjour3}       
\smartqed  
\usepackage{amsfonts}

\usepackage{latexsym}
\usepackage{amssymb}
\usepackage{amsmath}
\usepackage{graphicx}
\usepackage{hyperref}
\usepackage{caption}
\usepackage{subcaption}

\makeatletter
\def\@bibitem#1{\item[]%
    \if@filesw\immediate\write\@auxout{\string \bibcite {#1}{\the\value{\@listctr }}}\fi\ignorespaces}
\makeatletter

\newtheorem{Lemma}{Lemma}

\newtheorem{Theorem}[Lemma]{Theorem}

\renewcommand{\qed}{\hfill{\ \ \rule{2mm}{2mm}} \vspace{0.2in}}

\newcommand{\ind}{1\hspace{-2.3mm}{1}}

\begin{document}

\title{Euclidean traveling salesman problem with location dependent and power weighted edges}
\titlerunning{Traveling salesman problem}
\authorrunning{G. Ganesan}
\author{ \textbf{Ghurumuruhan Ganesan}
\thanks{E-Mail: \texttt{gganesan82@gmail.com} } \\
\ \\
Institute of Mathematical Sciences, HBNI, Chennai.}
\date{}
\maketitle

\begin{abstract}
Consider~\(n\) nodes~\(\{X_i\}_{1 \leq i \leq n}\)
independently distributed in the unit square~\(S,\)
each according to a distribution~\(f\)
and let~\(K_n\) be the complete graph formed by joining
each pair of nodes by a straight line segment.
For every edge~\(e\) in~\(K_n\) we associate a weight~\(w(e)\)
that may depend on the \emph{individual locations} of the endvertices
of~\(e\) and is not necessarily a power of the Euclidean length of~\(e.\)
Denoting~\(TSP_n\) to be the minimum weight of a spanning cycle of~\(K_n\)
corresponding to the travelling salesman problem (TSP)
and assuming an equivalence condition on the weight function~\(w(.),\)
we prove that~\(TSP_n\) appropriately
scaled and centred converges to zero a.s.\
and in mean as~\(n \rightarrow \infty.\)
We also obtain upper and lower bound deviation estimates for~\(TSP_n.\)




\vspace{0.1in} \noindent \textbf{Key words:} Traveling salesman problem, location dependent edge weights, deviation estimates.

\vspace{0.1in} \noindent \textbf{AMS 2000 Subject Classification:} Primary:
60J10, 60K35; Secondary: 60C05, 62E10, 90B15, 91D30.
\end{abstract}

\bigskip

\setcounter{equation}{0}
\renewcommand\theequation{\thesection.\arabic{equation}}
\section{Introduction}\label{intro}
The Traveling Salesman Problem (TSP) is the study of finding the minimum weight cycle containing all
the nodes of a graph where each edge is assigned a certain weight.
Often, the location of the nodes by itself is random
and the edge weights are taken to be the Euclidean
distance between the nodes. Thus the edge weight is a metric and for such models,
there is extensive literature dealing with various properties of the TSP
including the variance, deviation and almost sure convergence.
Beardwood et al (1959) use subadditive techniques to determine a.s.\ convergence
of the TSP length appropriately scaled to a constant~\(\beta_{TSP}\)
and recently Steinberger (2015) has obtained improved bounds for~\(\beta_{TSP}.\)
We refer to the books by Gutin and Punnen~(2006), Steele~(1997) and Yukich~(1998)
and papers by Steele (1981), Rhee (1993) and references therein
for various other aspects of TSP like the variance and complete convergence.

In this paper, we assume that the edge weight function is not necessarily a metric
and in fact might depend on the location of the endvertices of the edge.
This arises for example, in broadcasting problems of wireless ad-hoc networks where
each node is constrained to broadcast any received message at most once.
The cost of transmitting a message from a node might depend on the geographical conditions
surrounding and the node and the goal is therefore to
minimize the total cost incurred in sending packets to all nodes of the network.

In the rest of this section, we briefly describe the model under consideration
and state our result Theorem~\ref{tsp_thm} regarding the deviation estimates for the TSP length.


\subsection*{\em Model Description}
Let~\(f\) be any distribution on the unit square~\(S\) satisfying such that
\begin{equation}\label{f_eq}
\epsilon_1 \leq \inf_{x \in S} f(x) \leq \sup_{x \in S} f(x) \leq \epsilon_2
\end{equation}
for some constants~\(0 < \epsilon_1 \leq \epsilon_2 < \infty.\)
Throughout all constants are independent of~\(n.\)



Let~\(\{X_i\}_{i \geq 1}\) be independently
and identically distributed (i.i.d.) with the distribution~\(f(.)\) defined on the probability space \((\Omega, {\cal F}, \mathbb{P}).\)
For~\(n \geq 1,\) let~\(K_n = K(X_1,\ldots,X_n)\) be the complete graph whose edges are
obtained by connecting each pair of nodes~\(X_i\) and~\(X_j\) by the straight line segment~\((X_i,X_j)\)
with~\(X_i\) and~\(X_j\) as endvertices.

A path~\({\cal P} = (Y_1,\ldots,Y_t)\) is a subgraph of~\(K_n\)
with vertex set~\(\{Y_{j}\}_{1 \leq j \leq t}\)
and edge set~\(\{(Y_{j},Y_{{j+1}})\}_{1 \leq j \leq t-1}.\)
The nodes~\(Y_1\) and~\(Y_t\) are said to be \emph{connected} by
edges of the path~\({\cal P}.\)

Let~\(Y_1,\ldots,Y_t \subset \{X_k\}_{1 \leq k \leq n}\) be~\(t\) distinct nodes.  The subgraph~\({\cal C}  = (Y_1,Y_2,\ldots,Y_t,Y_1)\) with
vertex set~\(\{Y_{j}\}_{1 \leq j \leq t}\)
and edge set~\(\{(Y_{j},Y_{{j+1}})\}_{1 \leq j \leq t-1} \cup \{(Y_t,Y_1)\}\)
is said to be a \emph{cycle}. If~\({\cal C}\) contains all the nodes~\(\{X_i\}_{1 \leq i \leq n},\)
then~\({\cal C}\) is said to be a \emph{spanning cycle} of~\(K_n.\)

In what follows we assign weights to edges of the graph~\(K_n\) and study minimum weight spanning cycles.

\subsection*{\em Travelling Salesman Problem}
For points~\(x,y \in S,\) we let~\(d(x,y)\) denote the Euclidean distance between~\(x\) and~\(y\) and let~\(h : S \times S \rightarrow (0,\infty)\) be a deterministic measurable function satisfying
\begin{equation}\label{eq_met1}
c_1 d(x,y) \leq h(x,y)  = h(y,x) \leq c_2 d(x,y) 
\end{equation}
for some positive constants~\(c_1,c_2.\) For~\(\alpha > 0\) a constant and for~\(1 \leq i < j \leq n\) we let~\(h^{\alpha}(X_i,X_j) = h^{\alpha}(e)\) denote the \emph{weight} of the edge~\(e = (X_i,X_j),\) with exponent~\(\alpha.\) We remark here that unlike the Euclidean distance function~\(d,\) the edge weight function~\(h\) is not necessarily a metric. Also throughout, the quantity~\(\alpha\) appears in the superscript of any term only as an exponent.

The weight of a cycle~\({\cal C}\) in the graph~\(K_n\) is defined to be the sum of the weights of the edges in~\({\cal C};\) i.e.,
\begin{equation}\label{len_cyc_def}
W({\cal C}) := \sum_{e \in {\cal C}} h^{\alpha}(e).
\end{equation}
Let~\({\cal C}_{n}\) be a spanning cycle of the graph~\(K_n\) satisfying
\begin{equation}\label{min_weight_cycle}
TSP_n = W({\cal C}_{n}) := \min_{{\cal C}} W({\cal C}),
\end{equation}
where the minimum is taken over all spanning cycles~\({\cal C}.\) We refer to~\({\cal C}_{n}\) as the \emph{travelling salesman problem} (TSP) cycle with corresponding weight~\(TSP_n.\) If there is more than one choice for~\({\cal C}_{n},\) we choose one according to a deterministic rule.

Let~\(\epsilon_1,\epsilon_2\) be as in~(\ref{f_eq}) and set~\(\delta = \delta(\alpha) = \epsilon_1\) if the edge weight exponent is~\(\alpha \leq 1\) and~\(\delta = \epsilon_2\) if~\(\alpha > 1.\) Recalling that~\(c_1\) and~\(c_2\) are the bounds for the edge weight function~\(h\) as in~(\ref{eq_met1}), we define for~\(A > 0\) the terms~\(C_1(A)=C_1(A,\epsilon_1,\epsilon_2,\alpha)\)
and~\(C_2(A) = C_2(A,\epsilon_1,\epsilon_2,\alpha)\) as
\begin{eqnarray}
C_1(A)  &:=& \frac{(c_1A)^{\alpha}}{A^2}(1-e^{-\epsilon_1 A^2})e^{-8\epsilon_2 A^2} \text{ and } \nonumber\\
C_2(A)  &:=& (2c_2A)^{\alpha}\left(1 + \frac{1}{A^2}\left(\mathbb{E}\tilde{T}^{\alpha}+ \mathbb{E}\hat{T}^{\alpha}\right)\right), \label{c12def}
\end{eqnarray}
where~\(\tilde{T}\) is a geometric random variable with success parameter\\\(p = 1-e^{-\delta A^2}\left(1+\delta A^2 + \frac{\delta^2 A^4}{2}\right);\)
i.e.,~\(\mathbb{P}(\tilde{T}= k) = (1-p)^{k-1}p\) for all integers~\(k \geq 1\) and~\(\hat{T}\) is a geometric random variable with success parameter~\(1-p.\) We have the following result.
\begin{Theorem}\label{tsp_thm} Let~\(\alpha > 0\) be the edge weight exponent. For every~\(A > 0\) and every integer~\(k \geq 1\) and all~\(n \geq n_0(A,k,\epsilon_1,\epsilon_2,\alpha,c_1,c_2)\) large,
\begin{equation}\label{tsp_low_bounds}
\mathbb{P}\left(TSP_n \geq C_1(A) n^{1-\frac{\alpha}{2}}\left(1-\frac{4\sqrt{A}}{n^{1/4}}\right)\right) \geq 1-e^{-n^{1/3}},
\end{equation}
\begin{equation}\label{tsp_up_bounds}
\mathbb{P}\left(TSP_n \leq C_2(A) n^{1-\frac{\alpha}{2}}\left(1+ \frac{2}{n^{1/16}}\right)\right) \geq 1-\frac{1}{n^{2k}}
\end{equation}
and so
\begin{equation}
C_1^{k}(A) \left(1-\frac{37k\sqrt{A}}{n^{1/4}}\right)\leq \mathbb{E}\left(\frac{TSP^{k}_n}{n^{k\left(1-\frac{\alpha}{2}\right)}}\right) \nonumber\\
\leq C^{k}_2(A)\left(1+ \frac{3k}{n^{1/16}}\right). \label{exp_tsp_bound}
\end{equation}
\end{Theorem}
For the lower deviation estimates we determine the probability of favourable configurations that ensure a large enough number of relatively long edges in the TSP.
To prove the upper deviation estimate, we tile the unit square into small subsquares and join nodes within these subsquares to form an overall spanning cycle. We then estimate the length of this cycle via stochastic domination by homogenous processes (see Section~\ref{pf_tsp_thm_up2}).

\subsection*{\em Remarks on Theorem~\ref{tsp_thm}}
From Theorem~\ref{tsp_thm}, we see that the weights of the TSP in the location dependent case, is of the same order~\(n^{1-\frac{\alpha}{2}}\) as in the location independent case (see Steele (1988)). Using~(\ref{exp_tsp_bound}) we get that the normalized TSP weight~\(\frac{\mathbb{E} TSP_n}{n^{1-\frac{\alpha}{2}}}\) in fact satisfies
\begin{equation}\label{gr}
c_1^{\alpha} \cdot \beta_{low}(\alpha) \leq \liminf_n \frac{\mathbb{E} TSP_n}{n^{1-\frac{\alpha}{2}}} \leq \limsup_n \frac{\mathbb{E} MST_n}{n^{1-\frac{\alpha}{2}}} \leq c_2^{\alpha}\cdot \beta_{up}(\alpha),
\end{equation}
where
\begin{equation}\label{beta_low_def}
\beta_{low}(\alpha) = \beta_{low}(\alpha, \epsilon_1,\epsilon_2) := \sup_{A > 0} \frac{A^{\alpha}}{A^2}(1-e^{-\epsilon_1 A^2})e^{-8\epsilon_2 A^2},
\end{equation}
\begin{equation}\label{beta_up_def}
\beta_{up}(\alpha) = \beta_{up}(\alpha,\epsilon_1,\epsilon_2) := \inf_{A > 0} (2A)^{\alpha}\left(1 + \frac{\mathbb{E}{T_a}^{\alpha} + \mathbb{E}T_b^{\alpha}}{A^2}\right),
\end{equation}
and~\(T_a\) is a geometric random variable with success parameter~\(p\) and~\(T_b\) is a geometric random variable with success parameter~\(1-p\) (see~(\ref{c12def})).

\begin{figure}[tbp]
\centering
\includegraphics[width=3.5in, clip=true]{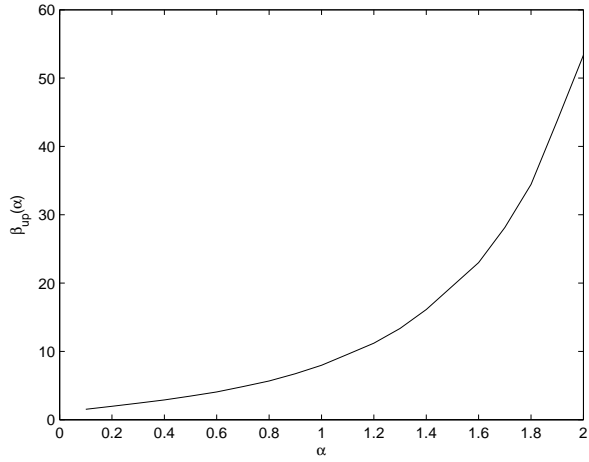}
\caption{Plot of~\(\beta_{up}(\alpha)\) as a function of~\(\alpha\) for the homogenous case~\(\epsilon_1 = \epsilon_2 = 1.\)}
\label{fig_up}
\end{figure}

\begin{figure}[tbp]
\centering
\includegraphics[width=3.5in, clip=true]{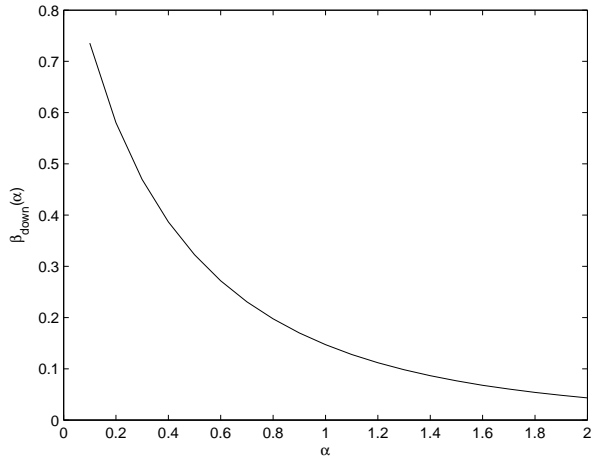}
\caption{Plot of~\(\beta_{down}(\alpha)\) as a function of~\(\alpha\) for the homogenous case~\(\epsilon_1 = \epsilon_2 = 1.\)}
\label{fig_down}
\end{figure}

For the case of homogenous distribution~\(\epsilon_1 = \epsilon_2 = 1\) we get
\begin{equation}\label{d_one}
\beta_{down}(\alpha) = \sup_{A > 0}A^{\alpha-2}(1-e^{-A^2})e^{-8A^2} > 0
\end{equation}
and
\begin{equation}\label{d_two}
\beta_{up}(\alpha) = \inf_{A > 0} (2A)^{\alpha} \left(1+\frac{\mathbb{E}T_a^{\alpha} + \mathbb{E}T_b^{\alpha}}{A^2}\right) < \infty,
\end{equation}
where~\(T_a\) is a geometric random variable with success parameter~\(p = 1-e^{-A^2}\left(1+A^2 + \frac{A^4}{2}\right)\)
and~\(T_b\) is a geometric random variable with success parameter~\(1-p.\) For illustration, we plot~\(\beta_{down}(\alpha)\) and~\(\beta_{up}(\alpha)\) as a function of~\(\alpha\) in Figures~\ref{fig_down} and~\ref{fig_up}, respectively. As we see from the figures~\(\beta_{up}(\alpha)\) increases with~\(\alpha\) and~\(\beta_{down}(\alpha)\) decreases with~\(\alpha.\)



As a final remark, we provide a simple upper bound for~\(\mathbb{E}T^{\alpha}\) for a geometric random variable~\(T\) with success parameter~\(p,\) in order to obtain quick evaluations of~\(\beta_{up}(\alpha)\) in~(\ref{beta_up_def}). Using~\(\mathbb{P}(T \geq k) = (1-p)^{k-1} \leq e^{-p(k-1)}\) we see that relation~(\ref{x_dist}) in Appendix is satisfied and so letting~\(r\) be the smallest integer greater than or equal~\(\alpha\) we have from~(\ref{disc_tel}) that
\[\mathbb{E}T^{\alpha} \leq \mathbb{E}T^{r} \leq \frac{r!}{(1-e^{-p})^{r}}.\] Plugging this estimate into~(\ref{beta_up_def}) provides
an upper bound for~\(\beta_{up}(\alpha).\)





The paper is organized as follows. In Section~\ref{pf_tsp_thm_up2},
we prove the deviation estimates in Theorem~\ref{tsp_thm}. In Section~\ref{pf_tsp_var_up} we obtain the variance upper bounds for~\(TSP_n\) and in Section~\ref{pf_tsp_var_low}, we obtain the variance lower bounds for~\(TSP_n.\) Combining the bounds we get that the variance grows roughly
of the order of~\(n^{1-\alpha}\) for the case~\(0 < \alpha < 1.\) Next in Section~\ref{pf_tsp_conv},
we prove the a.s.\ convergence for~\(TSP_n\) (Theorem~\ref{tsp_metric}), appropriately scaled and centred and finally
in Section~\ref{tsp_uni}, we obtain bounds for the scaled TSP weight when the nodes are uniformly distributed in the unit square (Theorem~\ref{tsp_unif}).





\setcounter{equation}{0}
\renewcommand\theequation{\thesection.\arabic{equation}}
\section{Proof of Theorem~\ref{tsp_thm}}\label{pf_tsp_thm_up2}
To prove the deviation estimates, we use Poissonization and let~\({\cal P}\) be a Poisson process
in the unit square~\(S\) with intensity~\(nf(.).\) We join  each pair
of nodes by a straight line segment and denote the resulting complete graph as~\(K^{(P)}_n.\)
We let~\(TSP_n^{(P)}\) be the TSP of~\(K^{(P)}_n\) by an analogous definition as in~(\ref{min_weight_cycle}). We first find deviation estimates for~\(TSP_n^{(P)}\) and  then use dePoissonization
to obtain corresponding estimates for~\(TSP_n,\) the TSP for the Binomial process as defined in~(\ref{min_weight_cycle}).

\begin{figure}[tbp]
\centering
\includegraphics[width=3in, trim= 20 200 50 110, clip=true]{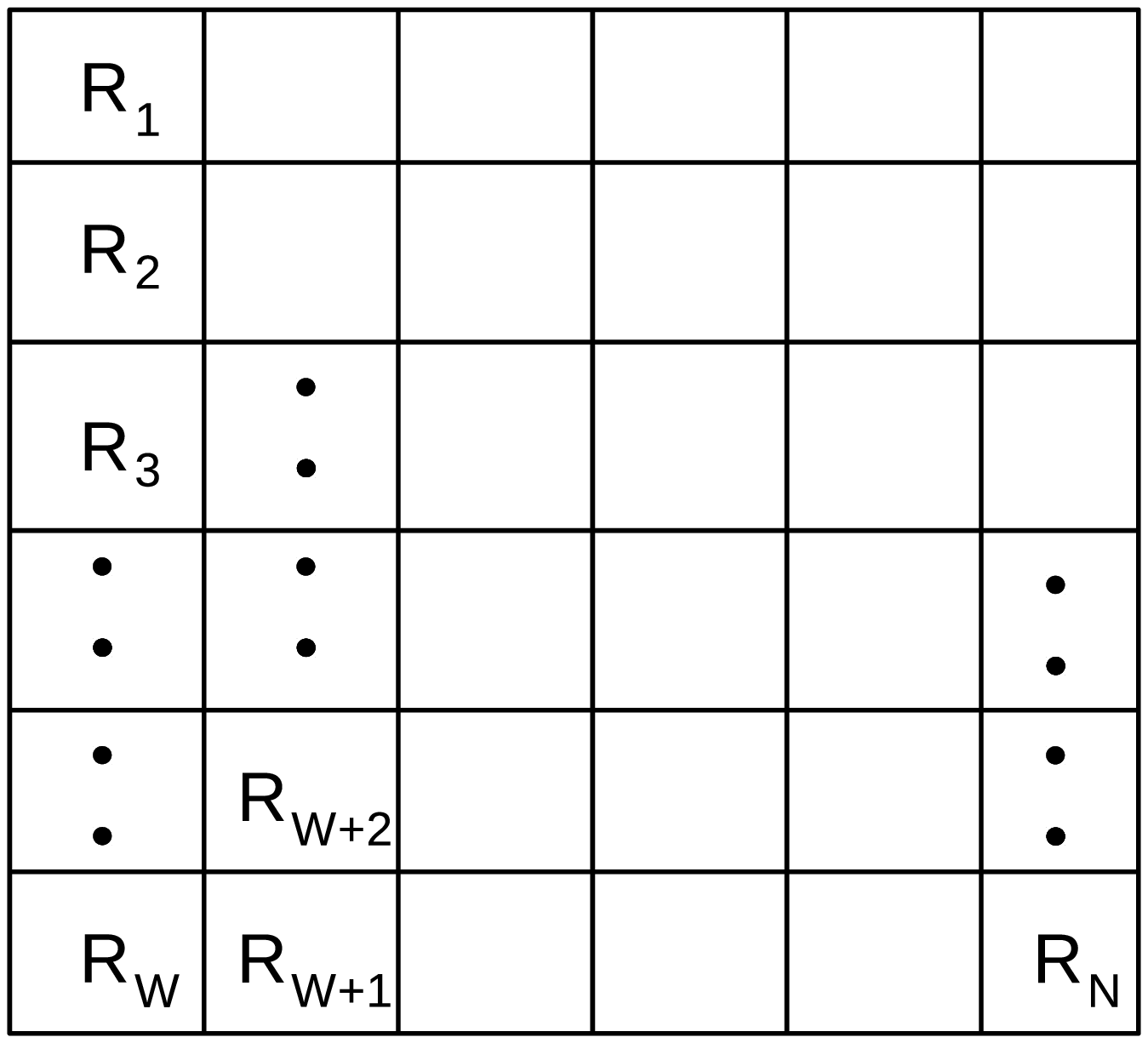}
\caption{Tiling the unit square into~\(N = \frac{n}{A^2}\) smaller~\(\frac{A}{\sqrt{n}} \times \frac{A}{\sqrt{n}}\) squares~\(\{R_l\}_{1 \leq l \leq \frac{n}{A^2}}.\)}
\label{fig_squares}
\end{figure}

For a  real number~\(A > 0,\) we tile the unit square~\(S\) into small~\(\frac{A(n)}{\sqrt{n}} \times \frac{A(n)}{\sqrt{n}}\) squares~\(\{R_i\}_{1 \leq i \leq \frac{n}{A^2(n)}}\) where~\(A(n) \in \left[A, A + \frac{1}{\log{n}}\right]\) is chosen such that~\(\frac{\sqrt{n}}{A(n)}\) is an integer. This is possible since
\begin{equation}
\frac{\sqrt{n}}{A} - \frac{\sqrt{n}}{A + (\log{n})^{-1}} = \frac{\sqrt{n}}{\log{n}}\cdot \frac{1}{A(A+(\log{n})^{-1})} \geq \frac{\sqrt{n}}{2A^2\log{n}} \label{poss}
\end{equation}
for all~\(n\) large. For notational simplicity, we denote~\(A(n)\) as~\(A\) henceforth and label the squares as in Figure~\ref{fig_squares} so that~\(R_i\) and~\(R_{i+1}\) share an edge for~\(1 \leq i \leq \frac{n}{A^2}-1.\)

\subsection*{\em Lower deviation bounds}
For~\(1 \leq i \leq \frac{n}{A^2}\) let~\(E(R_i)\) denote the event that the~\(\frac{A}{\sqrt{n}} \times \frac{A}{\sqrt{n}}\) square~\(R_i\) is occupied i.e., contains at least one node of~\({\cal P},\) and all squares sharing a corner with~\(R_i\) are empty. If~\(E(R_i)\) occurs, then there are at least two edges in the TSP of~\(K^{(P)}_n\) with one endvertex in~\(R_i\) and other endvertex in a square not sharing a corner with~\(R_i.\) Each such edge has a Euclidean length of at least~\(\frac{A}{\sqrt{n}}\) and so a weight of at least~\(\left(\frac{c_1A}{\sqrt{n}}\right)^{\alpha}\) (see~(\ref{eq_met1})). Consequently
\begin{equation}\label{low_bound}
TSP^{(P)}_n \geq \frac{1}{2} \cdot \sum_{i=1}^{\frac{n}{A^2}} 2\left(\frac{c_1A}{\sqrt{n}}\right)^{\alpha} \ind(E(R_i)) = \left(\frac{c_1A}{\sqrt{n}}\right)^{\alpha} \cdot G_{\alpha},
\end{equation}
where~\(G_{\alpha} := \sum_{i=1}^{\frac{n}{A^2}} \ind(E(R_i))\) and the factor~\(\frac{1}{2}\) occurs, since each edge is counted twice in the summation.

To estimate~\(G_{\alpha},\) we would like to split it into sums of independent r.v.s using the following construction. For a square~\(R_i,\) let~\({\cal N}(R_i)\) be the set of all squares sharing a corner with~\(R_i,\) including~\(R_i.\) If~\(R_i\) does not intersect the sides of the unit square~\(S,\) then there are~\(9\) squares in~\({\cal N}(R_i)\) and if~\(R_j\) is another square such that~\({\cal N}(R_i) \cap {\cal N}(R_j) = \emptyset,\) then the corresponding events~\(E(R_i)\) and~\(E(R_j)\) are independent, by Poisson property. We therefore extract nine disjoint subsets~\(\{{\cal U}_l\}_{1 \leq l \leq 9}\) of~\(\{R_i\}\) with the following properties:\\
\((A)\) If~\(R_i,R_j \in {\cal U}_l,\) then~\(\#{\cal N}(R_i) = \#{\cal N}(R_j) = 9\) and~\({\cal N}(R_i) \cap {\cal N}(R_j) = \emptyset.\)\\
\((B)\) The number of squares~\(\#{\cal U}_l \geq \frac{n}{9A^2} - \frac{4\sqrt{n}}{A}\) for each~\(1 \leq l \leq 9.\)\\
This is possible since there are at most~\(\frac{4\sqrt{n}}{A} - 4 < \frac{4\sqrt{n}}{A}\) squares in~\(\{R_k\}\) intersecting the sides of the unit square~\(S\) and the total number of squares in~\(\{R_k\}\) is~\(\frac{n}{A^2}.\)

We now write~\(G_{\alpha} = \sum_{i=1}^{\frac{n}{A^2}} \ind(E(R_i)) \geq \sum_{l=1}^{9} \sum_{R_i \in {\cal U}_l} \ind(E(R_i)),\) where each inner summation on the right side is a sum of independent Bernoulli random variables, which we bound via standard deviation estimates. Indeed for~\(1 \leq l \leq 9\) and~\(R_i \in {\cal U}_l,\) the number of nodes~\(N(R_i)\) is Poisson distributed with mean~\(n\int_{R_i}f(x)dx \in [\epsilon_1 A^2, \epsilon_2 A^2]\) (see~(\ref{f_eq})) and so~\(R_i\) is occupied with probability at least~\(1-e^{-\epsilon_1 A^2}.\) Also each of the eight squares sharing a corner with~\(R_i\) is empty with probability at least~\(e^{-\epsilon_2 A^2},\) implying that~\(\mathbb{P}(E(R_i)) \geq (1-e^{-\epsilon_1 A^2})e^{-8\epsilon_2 A^2}.\) Using the deviation estimate~(\ref{std_dev_down}) in Appendix with~\(\mu_1 = (1-e^{-\epsilon_1 A^2})e^{-8\epsilon_2 A^2}, m = \frac{n}{9A^2}-\frac{4\sqrt{n}}{A}\) and~\(\epsilon = \frac{1}{m^{1/4}}\) we then get that
\begin{equation}\label{ul_est}
\mathbb{P}_0\left(\sum_{R_i \in {\cal U}_l} \ind(E(R_i)) \geq (1-\epsilon)\left(\frac{n}{9A^2}-\frac{4\sqrt{n}}{A}\right)(1-e^{-\epsilon_1 A^2})e^{-8\epsilon_2 A^2}\right) \geq 1-e^{-D_1 \epsilon^2 n}
\end{equation}
for some constant~\(D_1 > 0\) not depending on~\(l.\) Since~\(m^{1/4} < \left(\frac{n}{9A^2}\right)^{1/4},\) we get that~\(D_1\epsilon^2 n \geq 2D_2 \sqrt{n}\)  for some constant~\(D_2 > 0\) and since~\(m^{1/4} > \left(\frac{n}{10A^2}\right)^{1/4}\) for all~\(n\) large, we have \[(1-\epsilon)\left(\frac{n}{9A^2} - \frac{4\sqrt{n}}{A}\right) \geq \frac{n}{9A^2} - \frac{4\sqrt{n}}{A} - \frac{n}{A^2 m^{1/4}} \geq \frac{n}{9A^2}\left(1 - \frac{36\sqrt{A}}{n^{1/4}}\right)\] for all~\(n\) large.

Letting
\[E_{low} := \left\{G_{\alpha}\geq (1-e^{-\epsilon_1 A^2})e^{-8\epsilon_2 A^2}\frac{n}{A^2}\left(1 - \frac{36\sqrt{A}}{n^{1/4}}\right)\right\},\] we get from~(\ref{ul_est}) that~\(\mathbb{P}_0(E_{low}) \geq 1-9e^{-2D_2 \sqrt{n}}\) and moreover, from~(\ref{low_bound}) we also get that \[TSP^{(P)}_n \ind(E_{low}) \geq \Delta_n := C_1(A) n^{1-\frac{\alpha}{2}} \left(1-\frac{36\sqrt{A}}{n^{1/4}}\right),\] where~\(C_1(A)\) is as defined in~(\ref{c12def}). From the estimate for the probability of the event~\(E_{low}\) above we therefore get~\(\mathbb{P}_0\left(TSP_n^{(P)}\geq \Delta_n\right) \geq 1- 9e^{-2D_2 \sqrt{n}}\) for all~\(n\) large. To convert the estimate from Poisson to the Binomial process, we let~\(E_P := \left\{TSP_n^{(P)} \geq \Delta_n\right\}, E := \left\{TSP_n \geq \Delta_n\right\}\) and use the dePoissonization formula ((Bradonjic et al. 2010), also proved below)
\begin{equation}\label{de_poiss_ax}
\mathbb{P}(E) \geq 1- D \sqrt{n} \mathbb{P}(E^c_P)
\end{equation}
for some constant~\(D > 0\) to get that~\(\mathbb{P}(E) \geq 1- D\sqrt{n}e^{-2D_2\sqrt{n}} \geq 1~-~e^{-D_2 \sqrt{n}}\) for all~\(n\) large. This proves~(\ref{tsp_low_bounds}) and so using~\[\mathbb{E}TSP^{k}_n \geq \mathbb{E}TSP^{k}_n \ind\left(MST_n \geq \Delta_n\right) \geq \Delta^{k}_n \left(1-e^{-D_2\sqrt{n}}\right)\] and
\[\left(1-\frac{36\sqrt{A}}{n^{1/4}}\right)^{k} (1-e^{-D_2 \sqrt{n}}) \geq \left(1-\frac{36k\sqrt{A}}{n^{1/4}}\right)(1-e^{-D_2 \sqrt{n}})  \geq 1-\frac{37k\sqrt{A}}{n^{1/4}}\] for all~\(n\) large, we also obtain the lower bound on the expectation in~(\ref{exp_tsp_bound}).

To see that~(\ref{de_poiss_ax}) is true, we let~\(N_P\) denote the random number of nodes of~\({\cal P}\) in all the squares~\(\{S_j\}\) so that~\(\mathbb{E}_0 N_P = n\) and~\(\mathbb{P}_0(N_P=n) = e^{-n}\frac{n^{n}}{n!} \geq \frac{D_1}{\sqrt{n}}\) for
some constant~\(D_1 > 0,\) using the Stirling formula. Given~\(N_P = n,\) the nodes of~\({\cal P}\)
are i.i.d.\ with distribution~\(f(.)\) as defined in~(\ref{f_eq}); i.e., \(\mathbb{P}_0(E_P^c|N_P = n)  = \mathbb{P}(E^c)\) and so
\[\mathbb{P}_0(E_P^c) \geq \mathbb{P}_0(E_P^c|N_P = n) \mathbb{P}_0(N_P = n) =
\mathbb{P}(E^c) \mathbb{P}_0(N_P = n) \geq \mathbb{P}(E^c)\frac{D_1}{\sqrt{n}},\] proving~(\ref{de_poiss_ax}).~\(\qed\)


\subsection*{\em Upper deviation bounds}
As before, we consider the Poisson process~\({\cal P}\) and first upper bound~\(TSP^{(P)}_n.\)  Intuitively, we would like to first connect all the nodes within each square~\(R_j\) to get a cycle and then merge all these cycles together to get an overall spanning cycle. However, it might so happen that some of the squares in~\(\{R_i\}\) may contain only one or two nodes and therefore, we take a slightly different approach. First we consider squares containing at least three nodes and join all these nodes to form a path. We then consider all the squares containing one or two nodes and connect these nodes together by a separate path and finally merge the two paths to get a spanning cycle.

Say that a square~\(R_i\) is a \emph{dense square} if it contains at least three nodes of~\({\cal P}\) and \emph{sparse} otherwise
and let~\(E_{bas}\) be the event that there are at least two sparse and at least two dense squares. To estimate the probability of  the event~\(E_{bas}\) we use the fact that the number of nodes~\(N(R_i)\) in the square~\(R_i\) is Poisson distributed with mean~\(\mathbb{E}_0N(R_i) = n \int_{R_i} f(x)dx \in [\epsilon_1 A^2, \epsilon_2 A^2]\) (see bounds for~\(f(.)\) in~(\ref{f_eq})). Therefore there are constants~\(p_0\) and~\(p_1\) not depending on~\(i\) such that
\begin{equation}\label{prob_spar}
0 < p_0 \leq p_{low} \leq p_{up} \leq p_1 < 1,
\end{equation}
where~\(p_{low} = \min\left(\mathbb{P}_0(R_i \text{ is dense}), \mathbb{P}_0(R_i \text{ is sparse}) \right)\)\\
and~\(p_{up} = \max\left(\mathbb{P}_0(R_i \text{ is dense}), \mathbb{P}_0(R_i \text{ is sparse}) \right).\)
Since the Poisson process is independent on disjoint sets,
there is at most one dense square among the~\(\frac{n}{A^2}\) squares in~\(\{R_i\}\) with probability
at most~\[(1-p_{low})^{\frac{n}{A^2}} + \frac{n}{A^2}(1-p_{low})^{\frac{n}{A^2}-1} \leq (1-p_{low})^{\frac{n}{4}-1}(n+1),\]
using~\(1 \leq A \leq 1+\frac{1}{\log{n}} \leq 2.\)
The same estimate holds for the event that there is at most one sparse square and so
we have that~
\begin{equation}\label{e_bas_est}
\mathbb{P}_0(E_{bas}) \geq 1-e^{-q_0n}
\end{equation}
for some constant~\(q_0 =q_0(\epsilon_1,\epsilon_2,A) > 0.\)

Henceforth we assume that~\(E_{bas}\) occurs and let~\(R_{i_1},R_{i_2},\ldots,R_{i_Q}, 1 \leq i_1 < i_2 <\ldots < i_Q \leq \frac{n}{A^2}, 2 \leq Q \leq \frac{n}{A^2}\) be all the dense squares. For~\(1 \leq j \leq Q,\) let~\({\cal P}_{i_j}\) be any spanning path containing all the nodes of~\(R_{i_j}\) and for~\(1 \leq j \leq Q-1\) let~\(e_{j+1}\) be any edge joining some node in~\(R_{i_j}\) and some node in~\(R_{i_{j+1}}\) so that the union~\({\cal P}_{dense} := \cup_{1 \leq j \leq Q} {\cal P}_{i_j} \cup \cup_{2 \leq l \leq Q} \{e_l\}\) is a path containing all the nodes in the dense squares.

For~\(1 \leq j \leq Q,\) there are~\(N(R_{i_j})\) nodes of the Poisson process in the square~\(R_{i_j}\) and any two such nodes are connected by an edge of Euclidean length at most~\(\frac{A\sqrt{2}}{\sqrt{n}}.\) Therefore the spanning path~\({\cal P}_{i_j}\) has a total weight of at most~\(N(R_{i_j}) \cdot \left(\frac{c_2A\sqrt{2}}{\sqrt{n}}\right)^{\alpha},\) using~(\ref{eq_met1}). The edge~\(e_{j+1}\) that connects some node in~\(R_{i_j}\) with some node of~\(R_{i_{j+1}}\) has a Euclidean length of at most~\(\frac{2T_{j+1}A}{\sqrt{n}}\) where~\(T_{j+1} := i_{j+1} - i_{j}\) and therefore has a weight of at most~\(\left(\frac{2c_2T_{j+1}A}{\sqrt{n}}\right)^{\alpha},\) again using~(\ref{eq_met1}). Setting~\(T_1 := i_1-1,T_{Q+1} := \frac{n}{A^2}-i_Q\) and~\(S_{\alpha} := \sum_{j=1}^{Q+1} T_j^{\alpha}\) we then get
\begin{eqnarray}
W\left({\cal P}_{dense}\right) &\leq& \sum_{j=1}^{Q} N(R_{i_j}) \cdot \left(\frac{c_2A\sqrt{2}}{\sqrt{n}}\right)^{\alpha} +
\sum_{j=2}^{Q} \left(\frac{2c_2T_{j}A}{\sqrt{n}}\right)^{\alpha} \nonumber\\
&\leq& \left(\frac{2c_2A}{\sqrt{n}}\right)^{\alpha}\left(\sum_{j=1}^{Q} N(R_{i_j}) + S_{\alpha}\right). \label{wl_dense2}
\end{eqnarray}

Suppose now that~\(R_{m_1},\ldots,R_{m_L}, 1 \leq m_1 < m_2 <\ldots < m_L \leq \frac{n}{A^2}, 2 \leq L \leq \frac{n}{A^2}\) are all the sparse squares. As before we connect the nodes in all these squares together to form a path~\({\cal P}_{sparse}\) whose weight is
\begin{equation}
W\left({\cal P}_{sparse}\right) \leq \sum_{l=1}^{L} N(R_{m_l}) \cdot \left(\frac{c_2A\sqrt{2}}{\sqrt{n}}\right)^{\alpha} +
\sum_{l=2}^{L} \left(\frac{2c_2U_{l}A}{\sqrt{n}}\right)^{\alpha}, \nonumber
\end{equation}
where~\(U_{l+1} := m_{l+1} - m_l\) for~\(1 \leq l \leq L-1.\) Setting~\(U_1 := m_1-1,U_{L+1} := \frac{n}{A^2}-m_L\) and~\(V_{\alpha} := \sum_{l=1}^{L+1} U_j^{\alpha}\) we then get
\begin{equation}
W\left({\cal P}_{sparse}\right) \leq \left(\frac{2c_2A}{\sqrt{n}}\right)^{\alpha}\left(\sum_{l=1}^{L} N(R_{m_l}) + V_{\alpha}\right). \label{wl_sparse}
\end{equation}
Adding~(\ref{wl_dense2}) and~(\ref{wl_sparse}) we get
\[W({\cal P}_{dense} ) + W({\cal P}_{sparse}) \leq \left(\frac{2c_2A}{\sqrt{n}}\right)^{\alpha}\left(\sum_{j=1}^{Q} N(R_{i_j}) + S_{\alpha} + \sum_{l=1}^{L} N(R_{m_l}) + V_{\alpha}\right)\]
and using the fact that the squares~\(\{R_{i_j}\} \cup \{R_{m_l}\}\) contain all the nodes of the Poisson process~\({\cal P},\) we then get
\begin{equation}
W({\cal P}_{dense} ) + W({\cal P}_{sparse}) \leq \left(\frac{2c_2A}{\sqrt{n}}\right)^{\alpha}\left(N_{tot} + S_{\alpha} + V_{\alpha}\right),\label{sum_length}
\end{equation}
where~\(N_{tot} := \sum_{i=1}^{\frac{n}{A^2}} N(R_{i})\) is the total number of nodes of~\({\cal P}\) in the unit square~\(S.\)

The path~\({\cal P}_{dense}\) has one endvertex~\(u_1\) in the square~\(R_{i_1}\) and the other endvertex~\(u_2\) in the square~\(R_{i_Q}.\) Similarly, the path~\({\cal P}_{sparse}\) has one endvertex~\(v_1\) in the square~\(R_{m_1}\) and the other endvertex~\(v_2\) in the square~\(R_{m_L}.\) Since the event~\(E_{bas}\) occurs, there are at least two dense and two sparse squares and so the vertices~\(u_1,u_2,v_1,v_2\) are all distinct. We join~\(u_1\) and~\(v_1\) by an edge~\(f_1\) and~\(u_2\) and~\(v_2\) by an edge~\(f_2\) so that the union~\({\cal C}_{tot} = {\cal P}_{dense} \cup {\cal P}_{sparse} \cup \{f_1,f_2\}\) is a spanning cycle containing all the nodes of the Poisson process~\({\cal P}.\)

Using the bounds for the edge weight function~\(h\) in~(\ref{eq_met1}), the weight of the edge~\(f_1\) is at most~\(\left(\frac{2c_2|i_1-m_1|A}{\sqrt{n}}\right)^{\alpha}\) and the weight of~\(f_2\) is at most~\(\left(\frac{2c_2|i_Q-m_L|A}{\sqrt{n}}\right)^{\alpha}.\) Therefore from~(\ref{sum_length}) we get that
\begin{eqnarray}
W({\cal C}_{tot}) &\leq& W({\cal P}_{dense}) + W({\cal P}_{sparse}) + \left(\frac{2c_2A}{\sqrt{n}}\right)^{\alpha} \cdot Z_{\alpha},\nonumber\\
&\leq& \left(\frac{2c_2A}{\sqrt{n}}\right)^{\alpha}\left(N_{tot}  + S_{\alpha} + V_{\alpha} + Z_{\alpha}\right),\label{up_bd1}
\end{eqnarray}
where~\(Z_{\alpha} := |i_1-m_1|^{\alpha} + |i_Q-m_L|^{\alpha}.\) From~(\ref{up_bd1}) and the fact that~\(TSP^{(P)}_n \leq W({\cal C}_{tot}),\) we get that
\begin{equation}
TSP^{(P)}_n \ind(E_{bas}) \leq \left(\frac{2c_2A}{\sqrt{n}}\right)^{\alpha}\left(N_{tot}  + S_{\alpha} + V_{\alpha}+Z_{\alpha}\right).\label{up_bdk}
\end{equation}

We evaluate each of the four terms in~(\ref{up_bdk}) separately. The first term~\(N_{tot}\) is a Poisson random variable with mean~\(n\) since this denotes the total number of nodes of the Poisson process in the unit square. From the deviation estimate~(\ref{std_dev_up}) in Appendix with~\(m=1,\mu_2 = n\) and~\(\epsilon = \frac{\log{n}}{\sqrt{n}},\) we have that~
\begin{equation}\label{e_node_est}
\mathbb{P}_0\left(N_{tot} \leq n\left(1+\frac{\log{n}}{\sqrt{n}}\right)\right) \geq 1-e^{-C(\log{n})^2}
\end{equation}
for some constant~\(C > 0.\) Setting~\(E_{node} :=  \left\{N_{tot} \leq n\left(1+\frac{\log{n}}{\sqrt{n}}\right)\right\},\) we get from~(\ref{up_bd1}) that
\begin{eqnarray}
&&TSP_n^{(P)}\ind(E_{bas} \cap E_{node})  \nonumber\\
&&\;\;\;\;\;\leq\;\;\left(\frac{2c_2A}{\sqrt{n}}\right)^{\alpha} \left(n^{}\left(1+\frac{\log{n}}{\sqrt{n}}\right)  + S_{\alpha} + V_{\alpha} + Z_{\alpha}\right).\label{up_bd2}
\end{eqnarray}

The following Lemma estimates~\(Z_{\alpha}, S_{\alpha}\) and~\(V_{\alpha}.\)
\begin{Lemma}\label{z_alp_lem}
Let~\(m \geq 2\) be any even integer constant. There exists a constant~\(D = D(m,\epsilon_1,\epsilon_2,\alpha) > 0\) such that
\begin{equation}\label{ez_est}
\mathbb{P}_0\left(Z_{\alpha} \leq D\left(\log{n}\right)^{\alpha}\right) \geq 1-\frac{4}{n^{m}},
\end{equation}
\begin{equation}
\mathbb{P}_{0}\left(S_{\alpha} \leq  \left(1+\frac{1}{n^{1/16}}\right)\frac{n}{A^2} \mathbb{E}{T}_a^{\alpha}\right) \geq 1- \frac{D}{n^{7m/16}} \label{til_est2}
\end{equation}
and
\begin{equation}
\mathbb{P}_{0}\left(V_{\alpha} \leq  \left(1+\frac{1}{n^{1/16}}\right)\frac{n}{A^2} \mathbb{E}{T}_b^{\alpha}\right) \geq 1- \frac{D}{n^{7m/16}}. \label{til_est_v}
\end{equation}
where~\({T}_a\) is a geometric random variable with success parameter~\[p_{dense} = 1-e^{-A^2\delta}\left(1+A^2\delta + \frac{A^4\delta^2}{2}\right)\] and~\({T}_b\) is a geometric random variable with success parameter~\[p_{sparse} = 1-p_{dense} = e^{-\delta A^2}\left(1+A^2\delta + \frac{A^4\delta^2}{2} \right).\]
\end{Lemma}
We prove Lemma~\ref{z_alp_lem} at the end. Assuming Lemma~\ref{z_alp_lem}, we complete the proof of Theorem~\ref{tsp_thm}.

Setting~\(E_Z = \{Z_{\alpha} \leq D\left(\log{n}\right)^{\alpha}\}\) we then get from~(\ref{up_bd2}) that
\begin{eqnarray}
&&TSP_n^{(P)} \ind(E_{bas} \cap E_{node} \cap E_Z)  \nonumber\\
&&\;\;\;\;\leq\;\; \left(\frac{2c_2A}{\sqrt{n}}\right)^{\alpha} \left(n\left(1+\frac{\log{n}}{\sqrt{n}}\right)  + D_0(\log{n})^{\alpha} +  S_{\alpha} + V_{\alpha}\right), \nonumber\\
&&\;\;\;\;\leq\;\; \left(\frac{2c_2A}{\sqrt{n}}\right)^{\alpha} \left(n\left(1+\frac{2\log{n}}{\sqrt{n}}\right)  + S_{\alpha} + V_{\alpha}\right)
\label{up_bdz}
\end{eqnarray}
for all~\(n\) large.

Letting~\(E_{SV}\) be the intersection of the events in the left hand sides of~(\ref{til_est2}) and~(\ref{til_est_v}), we have from the estimate~(\ref{up_bdz}) that
\begin{eqnarray}
&&TSP_n^{(P)} \ind(E_{bas} \cap E_{node} \cap E_Z \cap E_{SV})  \nonumber\\
&&\;\;\;\;\;\leq\;\; \left(\frac{2c_2A}{\sqrt{n}}\right)^{\alpha} \left(n\left(1+\frac{2\log{n}}{\sqrt{n}}\right)   + n\left(1+\frac{1}{n^{1/16}}\right)\frac{1}{A^2}(\mathbb{E}{T}_a^{\alpha} + \mathbb{E}{T}^{\alpha}_b)\right).\nonumber
\end{eqnarray}
Using~~\(\frac{2\log{n}}{\sqrt{n}} \leq \frac{1}{n^{1/16}}\) for all~\(n\) large, we then get
\begin{eqnarray}
TSP_n^{(P)} \ind(E_{bas} \cap E_{node} \cap E_Z \cap E_{SV})  \leq C_2(A) n^{1-\frac{\alpha}{2}}\left(1+ \frac{1}{n^{1/16}}\right), \label{up_bd3}
\end{eqnarray}
where~\(C_2(A)\) is as in~(\ref{c12def}). From~(\ref{up_bd3}) and the estimates for the events\\
\(E_{bas}, E_{node},E_{Z}\) and~\(E_{SV}\) from~(\ref{e_bas_est}),(\ref{e_node_est}) and~(\ref{ez_est}) and~(\ref{til_est2}), respectively, we have
\begin{eqnarray}
&&\mathbb{P}_0\left(TSP_n^{(P)} \leq C_2(A) n^{1-\frac{\alpha}{2}} \left(1+\frac{1}{n^{1/16}}\right)\right) \nonumber\\
&&\;\;\;\;\;\;\geq\;\;\;\; 1- e^{-q_0n} - e^{-D_5(\log{n})^2} - \frac{1}{n^{m}} - \frac{D_5}{n^{7m/16}} \nonumber\\
&&\;\;\;\;\;\;\geq\;\;\;\; 1-  \frac{2D_5}{n^{7m/16}} \label{main_poi_up}
\end{eqnarray}
for all~\(n\) large and some constant~\(D_5 > 0.\)

From~(\ref{main_poi_up}) and the dePoissonization formula~(\ref{de_poiss_ax}), we obtain
\begin{equation}
\mathbb{P}\left(TSP_n \leq C_2(A) n^{1-\frac{\alpha}{2}} \left(1+\frac{1}{n^{1/16}}\right)\right) \geq 1- \frac{D_6 \sqrt{n}}{n^{7m/16}} \label{temp1a}
\end{equation}
for some constant~\(D_6 > 0.\) Choosing~\(m\) large enough such that~\(\frac{7m}{16}-\frac{1}{2} \geq k,\) we obtain the estimate in~(\ref{tsp_up_bounds}). For bounding the expectation, we let~\(\Delta_n = C_2(A) n^{1-\frac{\alpha}{2}} \left(1+\frac{1}{n^{1/16}}\right)\) and write
~\(\mathbb{E}TSP_n^{k} = I_1 + I_2,\)
where
\begin{equation}\label{i1_est}
I_1 := \mathbb{E}TSP^{k}_n\ind(TSP_n \leq \Delta_n) \leq \Delta^{k}_n \leq C_2^{k}(A)n^{k\left(1-\frac{\alpha}{2}\right)}\left(1+\frac{2k}{n^{1/16}}\right),
\end{equation}
using~\(\left(1+\frac{1}{n^{1/16}}\right)^{k} \leq e^{k/n^{1/16}} \leq 1+\frac{2k}{n^{1/16}}\) for all~\(n\) large
and\\\(I_2 := \mathbb{E}TSP^k_n \ind(TSP_n > \Delta_n).\)

To evaluate~\(I_2,\) we use the estimate~\(TSP_n \leq n \cdot \left(c_2\sqrt{2}\right)^{\alpha}\) by~(\ref{eq_met1}), since there are~\(n\) edges in the spanning cycle and each such edge has an Euclidean length of at most~\(\sqrt{2}.\) Letting~\(\theta_m = \frac{7m}{16}-\frac{1}{2} - \frac{\alpha k}{2}\) and using the probability estimate~(\ref{temp1a})  we therefore get that
\begin{equation}\label{i2_est}
I_2  \leq \left(n \cdot \left(c_2\sqrt{2}\right)^{\alpha} \right)^{k} \cdot \frac{D_6 \sqrt{n}}{n^{7m/16}}  = n^{k\left(1-\frac{\alpha}{2}\right)}\frac{D_7}{n^{\theta_m}},
\end{equation}
for some constant~\(D_7 > 0.\) Adding~(\ref{i1_est}) and~(\ref{i2_est}) we get
\begin{equation}
\mathbb{E} TSP^{k}_n \leq C^{k}_2(A)n^{k\left(1-\frac{\alpha}{2}\right)}\left(1+\frac{2k}{n^{1/16}} + \frac{D_7}{n^{\theta_{m}}}\right) \nonumber
\end{equation}
and choosing~\(m\) larger if necessary so that~\(\theta_m \geq 1 > \frac{1}{16},\) we obtain the expectation upper bound in~(\ref{exp_tsp_bound}).~\(\qed\)

\emph{Proof of~(\ref{ez_est}) in Lemma~\ref{z_alp_lem}}: We write~\(Z_{\alpha} = |i_1-m_1|^{\alpha} + |i'_1-m'_1|^{\alpha},\) where~\(i'_1= \frac{n}{A^2}-i_Q\) and~\(m'_1 = \frac{n}{A^2}-m_L.\) Using~\((a+b)^{\alpha}\leq 2^{\alpha}(a^{\alpha} + b^{\alpha})\) for all~\(a,b,\alpha >0,\) we get~\(Z_{\alpha} \leq 2^{\alpha}(i_1^{\alpha} + m_1^{\alpha} + \left(i'_1\right)^{\alpha} + \left(m'_1\right)^{\alpha}).\) The term~\(i_1\) denotes the index of the first dense~\(\frac{A}{\sqrt{n}} \times \frac{A}{\sqrt{n}}\) square in~\(\{R_i\}.\) Since each square is independently dense with probability at least~\(p_0\) (see~(\ref{prob_spar})), we get that~\(\mathbb{P}_0(i_1 > l) \leq (1-p_0)^{l}.\) Setting~\(l = -\frac{m\log{n}}{\log(1-p_0)}\) where~\(m\) is an even integer constant to be determined later, we get~\(\mathbb{P}_0(i_1^{\alpha} > l^{\alpha}) \leq \frac{1}{n^{m}}.\) Analogous estimates hold for~\(m_1,i'_1\) and~\(m'_1\) and so we get~(\ref{ez_est}).~\(\qed\)

\emph{Proof of~(\ref{til_est2}) and~(\ref{til_est_v}) in Lemma~\ref{z_alp_lem}}: The quantities~\(S_{\alpha}\) and~\(V_{\alpha}\) are not i.i.d.\ sums and in what follows we evaluate~\(S_{\alpha}\) and an analogous analysis holds for~\(V_{\alpha}.\) The term~\(S_{\alpha}\) is defined for all configurations~\(\omega \in E_{bas}\) which contain at least two dense squares. We first extend the definition of~\(S_{\alpha}\) for all configurations~\(\omega\) as follows. If~\(\omega\) is such that there is exactly one dense square~\(R_{i_1}, 1 \leq i_1 \leq \frac{n}{A^2},\) then we set~\(S_{\alpha}(\omega) = (i_1-1)^{\alpha} + \left(\frac{n}{A^2}-i_1\right)^{\alpha}.\) If~\(\omega\) is such that there is no dense square, we set~\(S_{\alpha}(\omega) = \left(\frac{n}{A^2}-1\right)^{\alpha}.\) With this extended definition, we now show that~\(S_{\alpha}(\omega)\) is monotonic in~\(\omega\) in the sense that adding more nodes increases~\(S_{\alpha}\) if~\(\alpha \leq 1\) and decreases~\(S_{\alpha}\) if~\(\alpha > 1.\) This then allows us to use coupling and upper bound~\(S_{\alpha}\) by simply considering homogenous Poisson processes.

\underline{\emph{Monotonicity of~\(S_{\alpha}\)}}: Let~\(\omega\) be any configuration and suppose~\(\omega' = \omega \cup \{x\}\) is obtained by adding a single extra node at~\(x \in R_{j_0}\) for some~\(1 \leq j_0 \leq \frac{n}{A^2}.\) If~\(N(R_{j_0},\omega)\) and~\(N(R_{j_0},\omega')\) respectively denote the number of nodes of~\(\omega\) and~\(\omega'\) present in the square~\(R_{j_0},\) then~\(N(R_{j_0},\omega') = N(R_{j_0},\omega)+1.\)

If~\(N(R_{j_0},\omega) \leq 1,\) then~\(N(R_{j_0},\omega')\leq 2\) and so~\(R_{j_0}\) is not a dense square in the configuration~\(\omega'\) as well. This implies that~\(S_{\alpha}(\omega') = S_{\alpha}(\omega).\)  Analogous argument holds if~\(N(R_{j_0},\omega) \geq 3\) so that~\(R_{j_0}\) is already a dense square in the configuration~\(\omega.\) If on the other hand~\(N(R_{j_0},\omega) = 2,\) then~\(N(R_{j_0},\omega') = 3\) and so~\(R_{j_0}\) is a dense square in the configuration~\(\omega'\) but not in~\(\omega.\)

To compute~\(S_{\alpha}(\omega')\) in terms of~\(S_{\alpha}(\omega),\) suppose first that~\(\omega\) contains at least one dense square and let~\(1 \leq i_1(\omega) < \ldots < i_Q(\omega) \leq \frac{n}{A^2}\) be the indices of all the dense squares in~\(\{R_l\}.\) Setting~\(i_0(\omega) := 1\) and~\(i_{Q+1}(\omega) := \frac{n}{A^2},\) we then get that~\(S_{\alpha}(\omega) = \sum_{j=0}^{Q}(i_{j+1}(\omega)-i_{j}(\omega))^{\alpha}.\) There are only three possibilities for the index~\(j_0:\) Either~\(j_0=1\) or~\(j_0 = \frac{n}{A^2}\) or there exists~\(0 \leq a \leq Q\) such that~\(i_a(\omega)  < j_0 < i_{a+1}(\omega).\) If~\(j_0 \in \{1,\frac{n}{A^2}\},\) then again~\(S_{\alpha}(\omega) = S_{\alpha}(\omega').\) Else,~\(i_a(\omega) < j_0 < i_{a+1}(\omega)\) and so~
\begin{equation}\label{s_diff}
S_{\alpha}(\omega') = S_{\alpha}(\omega) + (i_{a+1}(\omega)-j_0)^{\alpha} +(j_0-i_a(\omega))^{\alpha} - (i_{a+1}(\omega)-i_a(\omega))^{\alpha}.
\end{equation}

If~\(\omega\) does not contain any dense square, then~\(S_{\alpha}(\omega) =  \left(\frac{n}{A^2}-1\right)^{\alpha}\)
and\\\(S_{\alpha}(\omega') = (j_0-1)^{\alpha} + \left(\frac{n}{A^2}-j_0\right)^{\alpha}\) so that~(\ref{s_diff}) is again satisfied with~\(i_a(\omega) = 1\) and~\(i_{a+1}(\omega)=\frac{n}{A^2}.\)  For~\(\alpha \leq 1\) we then use~\(y^{\alpha} + z^{\alpha} \geq (y+z)^{\alpha} \) for positive numbers~\(y,z\) in~(\ref{s_diff}), to get that~\(S_{\alpha}(\omega') \geq S_{\alpha}(\omega).\) If~\(\alpha > 1\) then~\(y^{\alpha}  + z^{\alpha} \leq (y+z)^{\alpha}\) and so~\(S_{\alpha}(\omega') \leq S_{\alpha}(\omega).\) This monotonicity property together with coupling allows us to upper bound~\(S_{\alpha}\) as follows. Letting~\(\delta = \epsilon_2\) if~\(\alpha \leq 1\) and~\(\delta = \epsilon_1\) if~\(\alpha > 1,\) we let~\({\cal P}_{\delta}\) be a homogenous Poisson process of intensity~\(\delta n\) on the unit square~\(S,\) defined on the probability space~\((\Omega_{\delta},{\cal F}_{\delta},\mathbb{P}_{\delta}).\)

Say that a square~\(R_i\) is a~\(\delta-\)dense square if~\(R_i\) contains at least three nodes of~\({\cal P}_{\delta}\) and let~\(F_{\delta}\) denote the event that there is at least one~\(\delta-\)dense square. As before we set~\(S^{(\delta)}_{\alpha} := \left(\frac{n}{A^2}-1\right)^{\alpha}\) if there is no~\(\delta-\)dense square. Suppose now that~\(F_{\delta}\) occurs and let~\(\{i^{(\delta)}_{j}\}_{1 \leq j \leq Q_{\delta}}\) be the indices of the~\(\delta-\)dense squares in~\(\{R_j\}.\) As in the definition of~\(S_{\alpha},\) we define~\(T^{(\delta)}_{j+1} := i^{(\delta)}_{j+1} - i^{(\delta)}_j\) for~\(1 \leq j \leq Q_{\delta}-1\) and set~\(T^{(\delta)}_1 := i^{(\delta)}_1-1\) and~\(T^{(\delta)}_{Q_{\delta}+1} := \frac{n}{A^2}-i^{(\delta)}_{Q_{\delta}}.\) Defining~\(S^{(\delta)}_{\alpha} := \sum_{j=1}^{Q_{\delta}+1} \left(T^{(\delta)}_j\right)^{\alpha},\) we have for any~\(x > 0\) that
\begin{equation}\label{mon_salpha}
\mathbb{P}_{\delta}\left(S^{(\delta)}_{\alpha} < x\right)\leq \mathbb{P}_0\left(S_{\alpha} < x\right).
\end{equation}
The proof follows by standard coupling techniques and for completeness, we have a proof in the Appendix.

To estimate~\(S^{(\delta)}_{\alpha}\) we let~\(N^{(\delta)}(R_i),1 \leq i \leq \frac{n}{A^2},\) be the random number of nodes of~\({\cal P}_{\delta}\) in the square~\(R_i.\) The random variables~\(\{N^{(\delta)}(R_i)\}\) are i.i.d.\ Poisson distributed each with mean~\(A^2\delta.\)  For~\(i \geq \frac{n}{A^2}+1,\) we define~\(N^{(\delta)}(R_i)\) to be i.i.d.\ Poisson random variables with mean~\(A^2\delta,\) that are also independent of~\(\{N^{(\delta)}(R_i)\}_{1 \leq i \leq \frac{n}{A^2}}.\) Without loss of generality, we associate the probability measure~\(\mathbb{P}_{\delta}\) for the random variables~\(\{N^{(\delta)}(R_i)\}_{i \geq \frac{n}{A^2} +1}\) as well.

Let~\(\tilde{T}_1 := \min\{j \geq 1 : N^{(\delta)}(R_j) \geq 3\}\)  and for~\(j \geq 2,\) let~\(\tilde{T}_j := \min\{k \geq \tilde{T}_{j-1}+1 : N^{(\delta)}(R_k) \geq 3\}-\tilde{T}_{j-1}.\) The random variables~\(\{\tilde{T}_i\}\) are nearly the same as~\(\{T^{(\delta)}_i\}\) in the following sense: We recall from discussion prior to~(\ref{mon_salpha}) that~\(F_{\delta}\) denotes the event that there exists at least one~\(\delta-\)dense square in~\(\{R_i\},\) containing at least three nodes of~\({\cal P}_{\delta}.\) If~\(F_{\delta}\) occurs, then~\(1 \leq Q_{\delta} \leq \frac{n}{A^2}\) and so~\(T^{(\delta)}_1 = \tilde{T}_1 -1, T^{(\delta)}_j = \tilde{T}_j\) for~\(2 \leq j \leq Q_{\delta}\) and~\(T^{(\delta)}_j \leq \tilde{T}_j\) for~\(j=Q_{\delta}+1.\) Consequently
\begin{equation}\label{dom1}
S_{\alpha}^{(\delta)}\ind(F_{\delta}) \leq \sum_{i=1}^{Q_{\delta}+1} \tilde{T}^{\alpha}_i\ind(F_{\delta}) \leq \sum_{i=1}^{\frac{n}{A^2}+1} \tilde{T}^{\alpha}_i
\end{equation}
since~\(Q_{\delta} \leq \frac{n}{A^2}.\) Also, arguing as in the proof of the estimate~(\ref{e_bas_est}) for the event~\(E_{bas}\) we have that~
\begin{equation}\label{f_del_est}
\mathbb{P}_{\delta}\left(F_{\delta}\right) \geq 1-e^{-qn}
\end{equation}
for some constant~\(q > 0.\)

The advantage of the above construction is that~\(\{\tilde{T}_i\}\) are i.i.d.\ geometric random variables with success parameter~\(p_{dense} = 1-e^{-A^2\delta}\left(1+A^2\delta + \frac{A^4\delta^2}{2}\right)\) and so all moments of~\(\tilde{T}_i^{\alpha}\) exist. Letting~\(\beta_i = \left(\tilde{T}_i^{\alpha} - \mathbb{E}_{\delta}\tilde{T}_i^{\alpha}\right)\) and~\(\beta_{tot} = \sum_{i=1}^{\frac{n}{A^2}+1} \beta_i \) we obtain for an even integer constant~\(m\) (to be determined later) that
\begin{equation}\label{e_del}
\mathbb{E}_{\delta}\beta_{tot}^m  = \mathbb{E}_{\delta} \sum_{(i_1,\ldots,i_m)} \beta_{i_1}\ldots\beta_{i_m}.
\end{equation}
For a tuple~\((i_1,\ldots,i_m)\) let~\(\{j_1,\ldots,j_w\}\) be the distinct integers in~\(\{i_1,\ldots,i_m\}\) with corresponding multiplicities~\(l_1,\ldots,l_w\) so that~\[\mathbb{E}_{\delta} \beta_{i_1}\ldots\beta_{i_m} = \mathbb{E}_{\delta} \beta_{j_1}^{l_1} \ldots \beta_{j_w}^{l_w} = \prod_{k=1}^{w} \mathbb{E}_{\delta} \beta_{j_k}^{l_k}.\] If~\(l_k = 1\) for some~\(1 \leq k \leq w,\) then~\(\mathbb{E}_{\delta} \beta_{i_1}\ldots\beta_{i_m} = 0\) and so for any non zero term in the summation in~(\ref{e_del}), there are at most~\(\frac{m}{2}\) distinct terms in~\(\{i_1,\ldots,i_m\}.\) This implies that
\[\mathbb{E}_{\delta}\beta_{tot}^{m} \leq D(m) {n \choose m/2} \leq D(m) n^{m/2}\] for some constant~\(D(m) > 0.\) For~\(\epsilon > 0\) we therefore get from Chebychev's inequality that
\[\mathbb{P}_{\delta}\left(|\beta_{tot}| > \epsilon \left(\frac{n}{A^2}+1\right) \mathbb{E}_{\delta}\tilde{T}_1^{\alpha}\right) \leq D_1\frac{\mathbb{E}_{\delta}(\beta_{tot}^m)}{n^m\epsilon^m} \leq \frac{D_2}{n^{m/2}\epsilon^m}\]
for some constants~\(D_1,D_2 > 0.\) Setting~\(\epsilon= \frac{1}{n^{1/16}}\) and using~\(\epsilon \left(\frac{n}{A^2}+1\right)  \leq \left(1+\frac{1}{n^{1/16}}\right) \frac{n}{A^2}\) for all~\(n\) large, we then get
\begin{equation}\label{til_est1}
\mathbb{P}_{\delta}\left(\sum_{i=1}^{\frac{n}{A^2}+1} \tilde{T}^{\alpha}_i \leq  \left(1+\frac{1}{n^{1/16}}\right)\frac{n}{A^2} \mathbb{E}\tilde{T}_1^{\alpha}\right) \geq 1-\frac{D_2}{n^{7m/16}}.
\end{equation}

From~(\ref{dom1}) and the estimate~(\ref{f_del_est}) for the event~\(F_{\delta},\) we further get
\begin{equation}
\mathbb{P}_{\delta}\left(S^{(\delta)}_{\alpha} \leq  \left(1+\frac{1}{n^{1/16}}\right)\frac{n}{A^2} \mathbb{E}\tilde{T}_1^{\alpha}\right) \geq 1-\frac{D_2}{n^{7m/16}}-e^{-q n} \geq 1- \frac{D_3}{n^{7m/16}} \nonumber
\end{equation}
for all~\(n\) large and some constant~\(D_3 > 0.\) Using the monotonicity relation~(\ref{mon_salpha}) we finally get~(\ref{til_est2}).
An analogous analysis yields~(\ref{til_est_v}).~\(\qed\)


\setcounter{equation}{0}
\renewcommand\theequation{\thesection.\arabic{equation}}
\section{Variance upper bound}\label{pf_tsp_var_up}
In this section, we obtain the variance upper bound for the overall TSP length~\(TSP_n\) as defined in~(\ref{min_weight_cycle}).
We have the following result.
\begin{Theorem}\label{tsp_var} Suppose the edge weight exponent is~\(0 < \alpha < 2.\) For every~\(\epsilon > \frac{(2-\alpha)(1+\alpha)}{2+\alpha},\) there is a constant~\(D_1 = D_1(\epsilon,\alpha,\epsilon_1,\epsilon_2) > 0\) such that~\[var(TSP_n) \leq D_1 n^{\epsilon}\] for all~\(n\) large.
\end{Theorem}
The desired bound is obtained via the martingale difference method that estimates the change in TSP lengths after adding or removing a single node.

We perform some preliminary computations.

\subsection*{\em One node difference estimates}
As in~(\ref{min_weight_cycle}), let~\(TSP_{n+1} = W({\cal C}_{n+1})\) be the weight of the TSP~\({\cal C}_{n+1}\) formed by the nodes~\(\{X_k\}_{1 \leq k \leq n+1}\) and for~\(1 \leq j \leq n+1,\) let~\(TSP_n(j) = W({\cal C}_n(j))\) be the weight of the TSP~\({\cal C}_{n}(j)\) formed by the nodes~\(\{X_k\}_{1 \leq k \neq j \leq n+1}.\) In this subsection, we find estimates for~\(|TSP_{n+1} - TSP_n(j)|,\) the change in the length of the TSP upon adding or removing a single node.

Suppose~\(v_j\) and~\(v_{j+1}\) are the neighbours of the node~\(X_j\) in the TSP~\({\cal C}_{n+1}\) so that the edges~\((v_j,X_j)\) and~\((X_j,v_{j+1})\) both belong to~\({\cal C}_{n+1}.\) Removing the edges~\((v_j,X_j)\) and~\((X_j,v_{j+1})\)
and adding the edge~\((v_j,v_{j+1}),\) we get a spanning cycle containing the nodes~\(\{X_k\}_{1 \leq k \neq j \leq n+1}\) and so
\begin{equation}\label{tsp_up22}
TSP_n(j) \leq TSP_{n+1} + h^{\alpha}(v_j,v_{j+1}) \leq TSP_{n+1} + c_2^{\alpha} d^{\alpha}(v_j,v_{j+1}),
\end{equation}
using~(\ref{eq_met1}).

By the triangle inequality we have that~\(d(v_j,v_{j+1}) \leq d(X_j,v_j) + d(X_j,v_{j+1})\) and so using~\((a+b)^{\alpha} \leq 2^{\alpha}(a^{\alpha}+b^{\alpha})\) for all~\(a,b,\alpha > 0\) we get that
\[d^{\alpha}(v_j,v_{j+1}) \leq 2^{\alpha}\left(d^{\alpha}(X_j,v_j) + d^{\alpha}(X_j,v_{j+1})\right) = 2^{\alpha} \sum_{v \in {\cal N}(X_j,{\cal C}_{n+1})} d^{\alpha}(X_j,v),\]
where~\({\cal N}(X_j,{\cal C}_{n+1})\) denotes the set of neighbours of~\(X_j\) in the spanning cycle~\({\cal C}_{n+1}.\) Thus
\begin{equation}\label{tsp_up2}
TSP_n(j) \leq TSP_{n+1} + g_1(X_j),
\end{equation}
where~\(g_1(X_j) := (2c_2)^{\alpha} \sum_{v \in {\cal N}(X_j,{\cal C}_{n+1})} d^{\alpha}(X_j,v).\)

To get an estimate in the reverse direction, we would like to identify edges of the TSP~\({\cal C}_n(j)\) with small length that are close to the (yet to be added) node~\(X_j.\) Recalling that~\(\epsilon > \frac{(2-\alpha)(1+\alpha)}{2+\alpha},\) we let~\(\theta >0\) be such that
\begin{equation}\label{theta_choice}
0 < \theta < \frac{\alpha}{2(2+\alpha)} \text{ and } 1-\epsilon < 2\theta\alpha < \frac{\alpha}{2}.
\end{equation}
To see such a~\(\theta\) exists, we use~\(\epsilon > \frac{(2-\alpha)(1+\alpha)}{2+\alpha}\) and choose~\(\theta > 0\) sufficiently close to but less than~\(\frac{\alpha}{2(2+\alpha)}\) so that~\(\epsilon > 1-2\theta\alpha > 1-2\left(\frac{\alpha}{2(2+\alpha)}\right)\alpha = \frac{(2-\alpha)(1+\alpha)}{2+\alpha}.\)
Fixing such a~\(\theta\) we use~\(\alpha < 2\) to get that~\(\theta < \frac{\alpha}{2(2+\alpha)} < \frac{1}{4}\) and therefore that~\(2\theta\alpha < \frac{\alpha}{2}.\)

Let~\(E_{good}(\theta)\) be the event that there is an edge~\((u,v) \in {\cal C}_n(j)\) such that each edge in the triangle formed by the nodes~\(X_j,u\) and~\(v\) has Euclidean length at most~\(\frac{1}{n^{\theta}}.\) If~\(E_{good}(\theta)\) occurs, then removing the edge~\((u,v)\) and adding the edges~\((X_j,u)\) and~\((X_j,v),\) we obtain a spanning cycle containing all the nodes~\(\{X_k\}_{1 \leq k \leq n+1}.\) The Euclidean length of each added edge is at most~\(\frac{1}{n^{\theta}}\) and so the weight of each added edge is at most~\(\left(\frac{c_2}{n^{\theta}}\right)^{\alpha},\) by~(\ref{eq_met1}).

If~\(E_{good}(\theta)\) does not occur, we then remove the shortest edge~\((u_1,v_1)\in {\cal C}_n(j)\) in terms of Euclidean length and add the edges~\((X_j,u_1)\) and~\((X_j,v_1).\) The Euclidean length of each of these newly added edges is at most~\(\sqrt{2}\) and so the corresponding weights are at most~\((c_2\sqrt{2})^{\alpha},\) each. Thus we get from the above discussion that~\(TSP_{n+1}\leq TSP_n(j) + g_2(X_j),\) where
\[g_2(X_j) := 2\left(\frac{c_2}{n^{\theta}}\right)^{\alpha} \ind(E_{good}(\theta)) + 2(c_2\sqrt{2})^{\alpha} \ind(E^{c}_{good}(\theta)).\] Combining with~(\ref{tsp_up2}), we get
\begin{equation}\label{cruc_tsp}
|TSP_{n+1} - TSP_n(j)| \leq g_1(X_j) + g_2(X_j).
\end{equation}

The following Lemma collects properties regarding~\(g_1\) and~\(g_2\) used later.\\
\begin{Lemma}\label{g_lem}
There is a constant~\(D > 0\) such that for all~\(n\) large,
\begin{equation}\label{g0_prop1}
\sum_{j=1}^{n+1} \mathbb{E}g_1(X_j) = (n+1)\mathbb{E}g_1(X_{n+1}) \leq Dn^{1-\frac{\alpha}{2}} ,
\end{equation}
\begin{equation}\label{g1_prop1}
g_1^2(X_j) \leq D\sum_{v \in {\cal N}(X_j,{\cal C}_{n+1})}d^{\alpha}(X_j,v), \sum_{j=1}^{n+1} \mathbb{E}g_1^2(X_j)  = (n+1)\mathbb{E}g_1^2(X_1) \leq 2Dn^{1-\frac{\alpha}{2}},
\end{equation}
\begin{equation}\label{g2_prop1}
\mathbb{P}(E_{good}(\theta)) \geq 1-\frac{1}{n^{2\theta\alpha}} \text{ and } \max_{1 \leq j \leq n+1}\mathbb{E}g_2^2(X_j) \leq \frac{D}{n^{2\theta \alpha}}.
\end{equation}
Moreover,
\begin{equation}\label{cruc_tsp_est}
\mathbb{E}|TSP_{n+1}-TSP_n| \leq \frac{2D}{n^{1-\epsilon}}.
\end{equation}
\end{Lemma}
The proof of~(\ref{cruc_tsp_est}) follows from~(\ref{cruc_tsp}) since
\[\mathbb{E}|TSP_{n+1}-TSP_n| \leq \mathbb{E}g_1(X_j) + \mathbb{E}g_2(X_j) \leq \frac{D}{n^{\frac{\alpha}{2}}} + \frac{D}{n^{2\theta\alpha}} \leq \frac{2D}{n^{1-\epsilon}},\]
since both~\(\frac{\alpha}{2}\) and~\(2\theta\alpha\) are at least~\(1-\epsilon,\) by our choice of~\(\theta\) in~(\ref{theta_choice}).\\\\
\emph{Proof of~(\ref{g0_prop1}) in Lemma~\ref{g_lem}}: Using~(\ref{eq_met1}), we have that
\[g_1(X_j) =  (2c_2)^{\alpha} \sum_{v \in {\cal N}(X_j,{\cal C}_{n+1})} d^{\alpha}(X_j,v) \leq \left(\frac{2c_2}{c_1}\right)^{\alpha}\sum_{v \in {\cal N}(X_j,{\cal C}_{n+1})} h^{\alpha}(X_j,v)\]
and so
\begin{equation}\label{fin_sum}
\sum_{j=1}^{n+1}\mathbb{E}g_1(X_j) \leq \left(\frac{2c_2}{c_1}\right)^{\alpha}\left(\sum_{j=1}^{n+1} \sum_{v \in {\cal N}(X_j,{\cal C}_{n+1})} h^{\alpha}(X_j,v)\right).
\end{equation}
The final double summation in~(\ref{fin_sum}) is simply twice the weight of the TSP cycle containing all the~\(n+1\) nodes and so using the expectation upper bound~(\ref{exp_tsp_bound}), we get the desired estimate in~(\ref{g0_prop1}).~\(\qed\)

\emph{Proof of~(\ref{g1_prop1}) in Lemma~\ref{g_lem}}: There are two edges containing~\(X_j\) as an endvertex in the cycle~\({\cal C}_{n+1}\) and so using~\((a+b)^2 \leq 2(a^2+b^2)\) we first get that \[g_1^2(X_j) \leq 2(2c_2)^{2\alpha} \sum_{v \in {\cal N}(X_j,{\cal C}_{n+1})} d^{2\alpha}(X_j,v).\] Since the Euclidean length of any edge is at most~\(\sqrt{2},\) we have that~\(d^{2\alpha}(X_j,v) \leq (\sqrt{2})^{\alpha}d^{\alpha}(X_j,v)\) and so
\[g_1^2(X_j) \leq 2(4c_2^2\sqrt{2})^{\alpha} \sum_{v \in {\cal N}(X_j,{\cal C}_{n+1})} d^{\alpha}(X_j,v),\] proving the first relation in~(\ref{g1_prop1}).
Consequently
\begin{eqnarray}
\sum_{j=1}^{n+1} \mathbb{E}g_1^2(X_j) &\leq& 2(4c_2^2\sqrt{2})^{\alpha} \mathbb{E}\sum_{j=1}^{n+1} \sum_{v \in {\cal N}(X_j,{\cal C}_{n+1})} d^{\alpha}(X_j,v) \nonumber\\
&\leq& 2\left(\frac{4c_2^2\sqrt{2}}{c_1}\right)^{\alpha} \mathbb{E}\sum_{j=1}^{n+1} \sum_{v \in {\cal N}(X_j,{\cal C}_{n+1})} h^{\alpha}(X_j,v),\label{g1_est1}
\end{eqnarray}
using~(\ref{eq_met1}). As before, the double summation in the right side of~(\ref{g1_est1}) is simply twice the weighted TSP length of~\({\cal C}_{n+1}\)
and so using expectation upper bound~(\ref{exp_tsp_bound}), we also get the second estimate in~(\ref{g1_prop1}).~\(\qed\)

\emph{Proof of~(\ref{g2_prop1}) in Lemma~\ref{g_lem}}: Let~\(a_0 < \frac{1}{16\sqrt{2}}\) be a positive constant and tile the unit square~\(S\) into~\(\frac{8A_0}{n^{\theta}} \times \frac{8A_0}{n^{\theta}}\) squares~\(Q_l, 1 \leq l \leq \frac{n^{2\theta}}{64A_0^2}\) where~\(A_0 = A_0(n) \in [a_0,2a_0)\) is such that~\(\frac{n^{\theta}}{8A_0}\) is an integer for all~\(n\) large.  If~\(Q^{0}_l\) is the~\(\frac{2A_0}{n^{\theta}} \times \frac{2A_0}{n^{\theta}}\) square with the same centre as~\(Q_l\) then the number of nodes~\(N^{0}_l\) of~\(\{X_i\}_{1 \leq i \neq j \leq n+1}\) within the square~\(Q_l^{0}\) is Binomially distributed with mean \[n\int_{Q^{0}_l} f(x)dx \geq \epsilon_1 n \cdot \frac{4A_0^2}{n^{2\theta}} \geq 4\epsilon_1 a_0^2 n^{1-2\theta},\] using the bounds for~\(f(.)\) in~(\ref{f_eq}).  From the deviation estimate~(\ref{std_dev_down}) in Appendix we then get that~\(N^{0}_l\) is at least~\(2\epsilon_1a_0^2 n^{1-2\theta}\) with probability at least\\\(1-e^{-2D_1n^{1-2\theta}},\) for some constant~\(D_1 > 0.\) Setting~\[E_{dense} := \bigcap_{1 \leq l \leq \frac{n^{2\theta}}{64A_0^2}} \{N^{0}_l \geq 2\epsilon_1 a_0^2 n^{1-2\theta}\},\] we therefore have
\begin{equation}\label{e_dense_est}
\mathbb{P}(E_{dense}) \geq 1-\frac{n^{2\theta}}{64A_0^2}e^{-2D_1n^{1-2\theta}} \geq 1-e^{-D_1n^{1-2\theta}}
\end{equation}
for all~\(n\) large.

The event~\(E_{dense}\) defined above is useful in the following way. Say that a \(\frac{8A_0}{n^{\theta}} \times \frac{8A_0}{n^{\theta}}\) square~\(Q_l\) is a \emph{bad} square if every edge of the TSP~\({\cal C}_n(j)\) with one endvertex within~\(Q^{0}_l,\) has its other endvertex outside~\(Q_l\) and let~\(E_{no\_bad}\) be the event that no square in~\(\{Q_k\}\) is bad. Suppose the event~\(E_{dense}\cap E_{no\_bad}\) occurs and the (new) node~\(X_j\) is added to the square~\(Q_l.\) There exists an edge~\(e = (u,v)\) of the TSP~\({\cal C}_n(j)\) with both endvertices in the smaller square~\(Q^0_l\) and so the Euclidean lengths~\(d(u,v),d(X_j,u)\) and~\(d(X_j,v)\) are each no more than~\(\frac{8A_0\sqrt{2}}{n^{\theta}},\) the length of the diagonal of~\(Q_l.\) Using~\(A_0 \leq 2a_0 < \frac{1}{8\sqrt{2}}\) (see the first sentence of this proof) we then that~\(\frac{8A_0\sqrt{2}}{n^{\theta}} \leq \frac{16a_0\sqrt{2}}{n^{\theta}} < \frac{1}{n^{\theta}}\) and so the event~\(E_{good}(\theta)\) occurs.


In what follows, we estimate the probability of the event~\(E_{dense} \cap E_{no\_bad}.\) Suppose that the event~\(E_{dense} \cap E^c_{no\_bad}\) occurs so that there is at least one bad square and let~\(Q_l\) be a bad square so that each node within the smaller subsquare~\(Q^0_l\) belongs to a (long) edge of the cycle~\({\cal C}_n(j)\) whose other endvertex lies outside the bigger square~\(Q_l.\) The Euclidean length of each such long edge is at least~\(\frac{3A_0}{n^{\theta}},\) the width of the annulus~\(Q_l \setminus Q^{0}_l\) and so using~\(A_0 \geq a_0,\) the weight of the long edge is at least~\(\left(\frac{3a_0c_1}{n^{\theta}}\right)^{\alpha},\) using the bounds for the weight function~\(h\) in~(\ref{eq_met1}). Additionally, since the event~\(E_{dense}\) occurs, there are at least~\(2\epsilon_1 a_0^2n^{1-2\theta}\) nodes in~\(Q^{0}_l\) and so
\begin{eqnarray}
TSP_n(j)\ind(E_{dense} \cap E^c_{no\_bad}) &\geq& 2\epsilon_1 a_0^2 n^{1-2\theta}\cdot \left(\frac{3a_0c_1}{n^{\theta}}\right)^{\alpha}\ind(E_{dense} \cap E^c_{no\_bad})  \nonumber\\
&=& 2\epsilon_1 a_0^2(3a_0c_1)^{\alpha} n^{1-(2+\alpha)\theta}\ind(E_{dense} \cap E^c_{no\_bad}).  \nonumber
\end{eqnarray}
For any integer constant~\(k \geq 1,\) we therefore get that
\begin{equation}
\mathbb{E}TSP^k_n(j)\ind(E_{dense} \cap E^c_{no\_bad}) \geq \left(2\epsilon_1a_0^2(3a_0c_1)^{\alpha}\right)^{k}n^{k(1-(2+\alpha)\theta)}\mathbb{P}(E_{dense} \cap E^c_{no\_bad}).  \nonumber
\end{equation}
But using the expectation upper bound in~(\ref{exp_tsp_bound}), we have
\[\mathbb{E}TSP^k_n(j)\ind(E_{dense} \cap E^c_{no\_bad}) \leq \mathbb{E}TSP^k_n(j) \leq D_1 n^{k\left(1-\frac{\alpha}{2}\right)}\] for some constant~\(D_1 >0\) and so~\(\mathbb{P}\left(E_{dense} \cap E^{c}_{no\_bad}\right) \leq  \frac{D_2}{n^{k\theta_1}},\)
where~\(D_2 > 0\) is a constant and~\(\theta_1 := \frac{\alpha}{2} - (2+\alpha)\theta.\) By choice~\(\theta < \frac{\alpha}{2(2+\alpha)}\) (see~(\ref{theta_choice})) and so choosing~\(k\) sufficiently large satisfying~\(k\theta_1 > 2\theta\alpha,\) we get~\[\mathbb{P}\left(E_{dense} \cap E^{c}_{no\_bad}\right) \leq  \frac{D_2}{n^{k\theta_1}} \leq \frac{1}{n^{2\theta\alpha}}\] for all~\(n\) large. Using the bounds for the probability of the event~\(E_{dense}\) from~(\ref{e_dense_est}), we then get that
\[\mathbb{P}(E_{dense} \cap E_{no\_bad}) \geq 1  - e^{-D_1n^{1-2\theta}}- \frac{1}{n^{2\theta\alpha}} \geq 1-\frac{2}{n^{2\theta\alpha}}\] for all~\(n\) large. This implies that~\(\mathbb{P}(E_{good}(\theta)) \geq 1-\frac{2}{n^{2\theta\alpha}}.\)~\(\qed\)



\subsection*{\em Variance upper bound for TSP}
We use one node difference estimate~(\ref{cruc_tsp}) together with the martingale difference method to obtain a bound for the variance. For~\(1 \leq j \leq n~+~1,\) let~\({\cal F}_j = \sigma\left(\{X_k\}_{1 \leq k \leq j}\right)\) denote the sigma field
generated by the node positions~\(\{X_k\}_{1 \leq k \leq j}.\) Defining the martingale difference
\begin{equation}\label{h_def}
H_j = \mathbb{E}(TSP_{n+1} | {\cal F}_j) - \mathbb{E}(TSP_{n+1} | {\cal F}_{j-1}),
\end{equation}
we then have that~\(TSP_{n+1} -\mathbb{E}TSP_{n+1} = \sum_{j=1}^{n+1} H_j\) and so by the martingale property
\begin{equation} \label{var_exp}
var(TSP_{n+1})  = \mathbb{E}\left(\sum_{j=1}^{n+1} H_j\right)^2 = \sum_{j=1}^{n+1} \mathbb{E}H_j^2.
\end{equation}

To evaluate~\(\mathbb{E}H_j^2\) we rewrite the martingale difference~\(H_j\) in a more convenient form.
Letting~\(X'_j\) be an independent copy of~\(X_j\) which is also independent of~\(\{X_k\}_{1 \leq k \neq j \leq n+1}\) we rewrite~
\begin{equation}\label{h_def2}
H_j = \mathbb{E}(TSP_{n+1}(X_j) - TSP_{n+1}(X'_j) | {\cal F}_j),
\end{equation}
where~\(TSP_{n+1}(X_j)\) is the weight of the MST formed by the nodes~\(\{X_i\}_{1 \leq i \leq n+1}\) and~\(TSP_{n+1}(X'_j)\)
is the weight of the TSP formed by the nodes\\\(\{X_i\}_{1 \leq i \neq j \leq n+1}~\cup~\{X'_j\}.\)
Using the triangle inequality and the one node difference estimate~(\ref{cruc_tsp}), we have
that~\(|TSP_{n+1}(X_j)-TSP_{n+1}(X'_j)|\) is bounded above as
\begin{eqnarray}
&&|TSP_{n+1}(X_j)-TSP_{n}(j)| + |TSP_{n+1}(X'_j)-TSP_{n}(j)| \nonumber\\
&&\;\;\;\leq\;\;g_1(X_j) + g_2(X_j) + g_1(X'_j) + g_2(X'_j), \nonumber
\end{eqnarray}
where~\(g_1\) and~\(g_2\) are as in~(\ref{cruc_tsp}) and we recall that~\(TSP_n(j)\) is the weight of the MST
formed by the nodes~\(\{X_k\}_{1 \leq k \neq j \leq n+1}.\)
Thus
\begin{eqnarray}
|H_j| &\leq& \mathbb{E}(|TSP_{n+1}(X_j)-TSP_{n+1}(X'_j)| | {\cal F}_j) \nonumber\\
&\leq& \mathbb{E}(g_1(X_j)|{\cal F}_j) + \mathbb{E}(g_2(X_j) |{\cal F}_j) + \mathbb{E}(g_1(X'_j)|{\cal F}_{j}) + \mathbb{E}(g_2(X'_j)|{\cal F}_{j}) \nonumber\\
&=& \mathbb{E}(g_1|{\cal F}_j) + \mathbb{E}(g_2 |{\cal F}_j) + \mathbb{E}(g_1|{\cal F}_{j-1}) + \mathbb{E}(g_2|{\cal F}_{j-1}). \nonumber
\end{eqnarray}
Using~\((a_1+a_2+a_3+a_4)^2 \leq 4(a_1^2 + a_2^2 + a_3^2 + a_4^2),\) we then get
\begin{eqnarray}
H_j^2 &\leq& 4\left(\left(\mathbb{E}(g_1|{\cal F}_j)\right)^2 + \left(\mathbb{E}(g_2 |{\cal F}_j)\right)^2 + \left(\mathbb{E}(g_1|{\cal F}_{j-1})\right)^2 + \left(\mathbb{E}(g_2|{\cal F}_{j-1})\right)^2\right)  \nonumber\\
&\leq& 4\left(\mathbb{E}(g^2_1|{\cal F}_j) + \mathbb{E}(g^2_2 |{\cal F}_j)+ \mathbb{E}(g^2_1|{\cal F}_{j-1}) + \mathbb{E}(g^2_2|{\cal F}_{j-1})\right) \nonumber
\end{eqnarray}
since~\((\mathbb{E}(X|{\cal F}))^2 \leq \mathbb{E}(X^2|{\cal F}).\) Thus~\(\mathbb{E}H_j^2 \leq 8\left(\mathbb{E}g_1^2(X_j) + \mathbb{E}g_2^2(X_j)\right)\) and plugging this in~(\ref{var_exp}), we have
\[var(TSP_{n+1}) = \sum_{j=1}^{n+1} \mathbb{E}H_j^2 \leq 8\left(\sum_{j=1}^{n+1}\mathbb{E}g_1^2(X_j) + \sum_{j=1}^{n+1} \mathbb{E}g_2^2(X_j)\right).\]
Plugging the estimates for~\(\mathbb{E}g_1^2(X_j)\) and~\(\mathbb{E}g_2^2(X_j)\) from~(\ref{g1_prop1}) and~(\ref{g2_prop1}), respectively,
we get~\(var(TSP_{n+1}) \leq D\left(n^{1-\frac{\alpha}{2}} + n^{1-2\theta\alpha}\right) \leq 2Dn^{\epsilon}\) for some constant~\(D > 0,\) by our choice of~\(\theta\) in~(\ref{theta_choice}).~\(\qed\)



\setcounter{equation}{0}
\renewcommand\theequation{\thesection.\arabic{equation}}
\section{Variance lower bound for~\(TSP_n\)}\label{pf_tsp_var_low}
\begin{Theorem}\label{var_tsp_low_thm}
Suppose there is a square~\(S_0\) with a constant side length~\(s_0\) such that~\(h(u,v) = d(u,v)\) if either~\(u\) or~\(v\) is in~\(S_0.\) For every~\(\alpha < 1,\) there is positive constant~\(D_2\) such that~\[var(TSP_n) \geq D_2n^{1-\alpha}\] for all~\(n\) large.
\end{Theorem}
We indirectly compute the variance lower bound for the TSP by computing a related spanning path~\({\cal P}_{app}(n+1)\) which has approximately the same length as the TSP
of the nodes~\(\{X_i\}_{1 \leq i \leq n+1}\) and estimate the variance of
the weight of~\({\cal P}_{app}(n+1).\) Using this estimate, we obtain the desired
variance lower bound for~\(TSP_{n+1}.\)

We tile the unit square~\(S\) into disjoint~\(\frac{400A}{\sqrt{n}} \times \frac{400A}{\sqrt{n}}\)
squares~\(\{R_i\}\) where~\(1 \leq A \leq 1+\frac{1}{\log{n}}\) is such that~\(\frac{\sqrt{n}}{400A}\)
is an integer for all~\(n\) large. This is possible by an argument following~(\ref{poss}).
Let~\(R(1) \in \{R_i\}\) be the (random) square containing the node~\(X_1\)
and let~\(R^{small}(1)\) be the~\(\frac{200A}{\sqrt{n}} \times \frac{200A}{\sqrt{n}}\) square with the same centre as~\(R(1)\)
and contained in~\(R(1).\)

Let~\({\cal P}_{in}(n+1)\) be the minimum weight spanning path formed by the nodes of~\(\{X_i\}_{1 \leq i \leq n+1}\)
present within the square~\(R^{small}(1)\) and let~\(u_{in}\) and~\(v_{in}\)
be the endvertices of~\({\cal P}_{in}(n+1).\) We define~\({\cal P}_{in}(n+1)\)
to be the minimum weight \emph{in-spanning path}. Let~\(v_{close}\)
be the node closest (in terms of Euclidean length) to but outside~\(R^{small}(1)\)
and suppose~\(v_{in}\) is closer in terms of Euclidean length to~\(v_{close}\) than~\(u_{in}.\)
We define~\(e_{min} = (v_{in},v_{close})\) to be the \emph{cross edge}. Finally, let~\({\cal P}_{out}(n+1)\) be
the minimum weight spanning path formed by the nodes of~\(\{X_i\}_{1 \leq i \leq n+1}\) present
outside~\(R^{small}(1)\) and having~\(v_{close}\) as an endvertex. We denote~\({\cal P}_{out}(n+1)\)
to be the minimum weight \emph{out-spanning path}. By construction, the union
\begin{equation}\label{p_app}
{\cal P}_{app}(n+1) = {\cal P}_{app}(n+1, \{X_i\}_{1 \leq i \leq n+1}) := {\cal P}_{in}(n+1) \cup \{e_{min}\} \cup {\cal P}_{out}(n+1)
\end{equation}
is an \emph{overall spanning path} containing all the nodes~\(\{X_i\}_{1 \leq i \leq n+1}\) and we define~\({\cal P}_{app}(n~+~1)\)
to be the \emph{approximate overall TSP path}.

As in~(\ref{len_cyc_def}), we define~\(W_{n+1} := W({\cal P}_{app}(n+1))=\sum_{e \in {\cal P}_{app}(n+1)} h^{\alpha}(e)\) to be the weight of the approximate overall TSP path~\({\cal P}_{app}(n+1)\) and recall that~\(TSP_{n+1}\) is the minimum weight of the overall spanning cycle containing all the~\(n+1\) nodes~\(\{X_i\}_{1 \leq i \leq n+1}.\) The following result shows that the approximate overall TSP path has nearly the same weight as the optimal overall TSP cycle and obtains a lower bound on the variance.
\begin{Lemma}\label{var_low_lem1} There exists a constant~\(D > 0\) such that
\begin{equation}\label{var_close}
\mathbb{E}|TSP_{n+1}-W_{n+1}|^2 \leq D(\log{n})^2
\end{equation}
and
\begin{equation}\label{var_p_app}
var(W_{n+1}) \geq D\cdot n^{1-\alpha}
\end{equation}
for all~\(n\) large.
\end{Lemma}
From the bounds in Lemma~\ref{var_low_lem1}, we obtain a lower bound for the variance of~\(TSP_{n+1}\) as follows. Using~\((a+b)^2 \leq 2(a^2+b^2)\) we have for any two random variables~\(X\) and~\(Y\) that~\[var(X+Y) \leq 2(var(X) + var(Y)) \leq 2var(X) + 2\mathbb{E}Y^2.\]
Setting~\(X = TSP_{n+1}, Y = W_{n+1} - TSP_{n+1}\) and using the bounds in Lemma~\ref{var_low_lem1}, we therefore get
\begin{eqnarray}
var(TSP_{n+1}) &\geq& \frac{1}{2}\cdot var(W_{n+1}) - \mathbb{E}|TSP_{n+1}-W_{n+1}|^2 \nonumber\\
&\geq& \frac{1}{2}\cdot D\cdot n^{1-\alpha} - D(\log{n})^2, \nonumber
\end{eqnarray}
obtaining the desired bound on the variance of~\(TSP_{n+1}.\)

\emph{Proof of~(\ref{var_close}) in Lemma~\ref{var_low_lem1}}:  Joining the endvertices of the approximate overall TSP path~\({\cal P}_{app}(n+1)\) by an edge, we get an overall spanning cycle containing all the nodes~\(\{X_i\}_{1 \leq i \leq n+1}.\) The added edge has a Euclidean length of at most~\(\sqrt{2}\) and so a weight of at most~\((c_2\sqrt{2})^{\alpha},\) by~(\ref{eq_met1}). Thus
\begin{equation}\label{x1_eq}
TSP_{n+1} \leq W_{n+1} + (c_2\sqrt{2})^{\alpha}.
\end{equation}

To get an estimate in the other direction, we first get from~(\ref{p_app}) that
\[W_{n+1} = W({\cal P}_{in}(n+1)) + h^{\alpha}(e_{min}) + W({\cal P}_{out}(n+1)) \] and estimate the lengths of~\({\cal P}_{in}(n+1), e_{min}\) and~\({\cal P}_{out}(n+1)\) in that order. If~\(N_{in}\) is the number nodes of~\(\{X_i\}_{1 \leq i \leq n+1}\) present in the square~\(R^{small}(1),\)
then the number of edges in the in-spanning path~\({\cal P}_{in}(n+1)\) is~\(N_{in}-1~<~N_{in}.\) As argued before, each such edge has a weight of at most~\((c_2\sqrt{2})^{\alpha}\) and so~\(W({\cal P}_{in}(n+1)) \leq N_{in}\cdot(c_2 \sqrt{2})^{\alpha}.\)
Similarly~\(h^{\alpha}(e_{min}) \leq (c_2\sqrt{2})^{\alpha}.\)

We now find bounds for~\(W({\cal P}_{out}(n+1))\) in terms of~\(TSP_{n+1}\) as follows. We recall that~\({\cal P}_{out}(n+1)\) is an out-spanning path with vertex set being the nodes of~\(\{X_i\}_{1 \leq i \leq n+1}\) outside~\(R^{small}(1).\) Also~\({\cal P}_{out}(n+1)\) has~\(v_{close}\) as an endvertex, where~\(v_{close}\) is the node closest to but outside~\(R^{small}(1).\) We therefore obtain an upper bound for~\(W({\cal P}_{out}(n+1))\) by constructing an out-spanning path~\({\cal P}_{fin}\) with~\(v_{close}\) as an endvertex, starting  from the overall spanning cycle~\({\cal C}_{n+1} = (u_1,\ldots,u_{n+1},u_1)\) with weight~\(TSP_{n+1}.\)

Let~\((u_{j_1},\ldots,u_{j_2}),(u_{j_3},\ldots,u_{j_4}),\ldots,(u_{j_s},\ldots,u_{j_{s+1}})\) be the subpaths of~\({\cal C}_{n+1}\) that are contained within the square~\(R^{small}(1)\) in the following sense: For example in the path~\((u_{j_1},\ldots,u_{j_2}),\) the nodes~\(u_{j_1}\) and~\(u_{j_2}\) are present outside~\(R^{small}(1)\) and the rest of the nodes are present inside~\(R^{small}(1).\) Removing all these subpaths and adding the edges~\((u_{j_1},u_{j_2}),(u_{j_3},u_{j_4}),\ldots,(u_{j_{s-1}},u_{j_s}),\) we therefore get an out-spanning cycle~\({\cal C}_{temp}.\) Further removing an edge containing~\(v_{close}\) from~\({\cal C}_{temp},\) we get an out-spanning path~\({\cal P}_{fin}\) having~\(v_{close}\) as an endvertex. Thus~\(W({\cal P}_{out}(n+1)) \leq W({\cal P}_{fin})\) and it suffices to upper bound~\(W({\cal P}_{fin}).\)


Recalling that~\(N_{in}\) is the number nodes of~\(\{X_i\}_{1 \leq i \leq n+1}\) present in~\(R^{small}(1),\) we get that the total number of edges added in the above procedure is at most~\(s \leq N_{in}.\) As argued before, the weight of any edge is at most~\((c_2\sqrt{2})^{\alpha}\) by~(\ref{eq_met1}) and so~\[W({\cal P}_{out}(n+1)) \leq W({\cal P}_{fin}) \leq TSP_{n+1} +  N_{in} \cdot (c_2\sqrt{2})^{\alpha}.\] Combining this with the estimates for~\(W({\cal P}_{in}(n+1))\)
and the weight of the edge~\(e_{min}\) obtained before,  we get
\begin{equation}\nonumber
W_{n+1} \leq TSP_{n+1}  + (2N_{in} + 1)\cdot (c_2\sqrt{2})^{\alpha}.
\end{equation}
and so from~(\ref{x1_eq}) we get
\begin{equation} \nonumber
|W_{n+1} - TSP_{n+1}| \leq (2N_{in} + 2)\cdot (c_2\sqrt{2})^{\alpha}.
\end{equation}
Taking expectations we get
\begin{equation}\nonumber
\mathbb{E}|W_{n+1} -TSP_{n+1}|^2 \leq 4\cdot (c_2\sqrt{2})^{2\alpha} \cdot \mathbb{E}(N_{in}+1)^2.
\end{equation}
and we prove below that~\(\mathbb{E}(N_{in}+1)^2 \leq 2(\log{n})^2,\) thus obtaining~(\ref{var_close}).

To estimate~\(N_{in},\) we write~\(N_{in} \leq \max_{i} N(R_i),\) where~\(N(R_i)\) is the number of nodes of~\(\{X_l\}_{1 \leq l \leq n+1}\) in the square~\(R_i.\) Any square in~\(\{R_i\}\) has an area of~\(\frac{10^{4}A^2}{n}\) and so using~(\ref{f_eq}) and~\(A \leq 1+\frac{1}{\log{n}} \leq 2,\) we get that~\(N(R_i)\) is Binomially distributed with a mean and a variance of at most~\(\epsilon_2 (n+1) \frac{10^{4}A^2}{n} \leq D_1\)
for some constant~\(D_1 \geq 1.\) Therefore from the proof of the deviation estimate~(\ref{std_dev_up}) in Appendix,
we get for~\(s > 0\) that~\(\mathbb{E}e^{sN(R_i)} \leq e^{D_1(e^{s}-1)}\) and so by Chernoff bound,
~\[\mathbb{P}(N(R_i) \geq \log{n}) \leq e^{-s\log{n}}e^{D_1(e^{s}-1)} \leq \frac{D_2}{n^{4}}\] by setting~\(s = 4\)
and~\(D_2 = e^{D_1(e^{4}-1)}.\) Since there are at most~\(\frac{n}{10^{4}A^2}\) squares in~\(\{R_i\}\) we further get
\begin{equation}\label{nri_est}
\mathbb{P}(\max_{i} N(R_i) \geq \log{n}) \leq \frac{n}{10^{4}A^2}\cdot \frac{D_2}{n^4} \leq \frac{D_3}{n^3}
\end{equation}
for some constant~\(D_3 > 0\) and so
\begin{eqnarray}
\mathbb{E}(N_{in}+1)^2 &=& \mathbb{E}(N_{in}+1)^2\ind(\max_{i}N(R_i) \leq \log{n}) \nonumber\\
&&\;\;\;\;\;+\;\;\mathbb{E}(N_{in}+1)^2 \ind(\max_{i} N(R_i) > \log{n}) \nonumber\\
&\leq& (\log{n}+1)^2 + \mathbb{E}(N_{in}+1)^2 \ind(\max_{i} N(R_i) > \log{n}). \nonumber
\end{eqnarray}
Using~\(N_{in} \leq n+1\) and~(\ref{nri_est}), we get
\begin{eqnarray}
\mathbb{E}(N_{in}+1)^2 \ind(\max_{i} N(R_i) > \log{n}) &\leq& (n+1)^2\mathbb{P}(\max_{i}N(R_i) >\log{n}) \nonumber\\
&\leq& (n+1)^2\cdot \frac{D_3}{n^3} \nonumber\\
&\leq& 1 \nonumber
\end{eqnarray}
for all~\(n\) large. Thus~\(\mathbb{E}(N_{in}+1)^2 \leq 2(\log{n})^2,\) proving the desired estimate.~\(\qed\)

\emph{Proof of~(\ref{var_p_app}) in Lemma~\ref{var_low_lem1}} We perform some preliminary computations. As in the proof of the variance upper bound for~\(1 \leq j \leq n+1\) we let~\({\cal F}_j\)
denote the sigma field generated by the nodes~\(\{X_i\}_{1 \leq i \leq j}\) and define the martingale difference~\(G_j = \mathbb{E}(W_{n+1}|{\cal F}_j) - \mathbb{E}(W_{n+1}|{\cal F}_{j-1})\) so that
\begin{equation}\label{var_wl_one}
var(W_{n+1}) = \left(\sum_{j=1}^{n+1} \mathbb{E}G_j\right)^{2} = \sum_{j=1}^{n+1} \mathbb{E}G_j^2,
\end{equation}
by the martingale property. To evaluate~\(\mathbb{E}G_j^2\) we rewrite the martingale difference~\(G_j\) in a more convenient form as
\begin{equation}
G_j = \mathbb{E}(W_{n+1}(X_j) - W_{n+1}(X'_j)|{\cal F}_{j}),\label{exp_hj}
\end{equation}
where~\(X'_j\) is an independent copy of~\(X_j,\) also independent of~\(\{X_i\}_{1 \leq i \neq j \leq n+1}\) and~\(W_{n+1}(X_j)\) and~\(W_{n+1}(X'_j)\) are the weights of the approximate TSP paths formed by the nodes~\(\{X_i\}_{1 \leq i \leq n+1}\) and~\(\{X_i\}_{1 \leq i \neq j \leq n+1} \cup \{X'_j\}.\) In words,\\\(W_{n+1}(X_j) - W_{n+1}(X'_j)\) represents the change in the length of the approximate overall TSP path after ``replacing" the node~\(X_j\) by~\(X'_j.\)

Below we define an event~\(E_{good}(j)\) with the following properties:\\
\((p1)\) There exists a constant~\(D >0\) not depending on~\(j\) such that~\[(W_{n+1}(X_j)-W_{n+1}(X'_j)) \ind(E_{good}(j)) \geq D \cdot n^{-\frac{\alpha}{2}}\cdot \ind(E_{good}(j)).\]
\((p2)\) The minimum probability~\(\min_{2 \leq j \leq n+1} \mathbb{P}(E_{good}(j)) \geq c_0\) for some constant~\(c_0  >0.\)\\
Assuming~\((p1)-(p2),\) we get from the expression for~\(G_j\) in~(\ref{exp_hj}) that
\[\mathbb{E}|G_j| \geq \mathbb{E}|G_j| \ind(E_{good}(j)) = \mathbb{E}G_j \ind(E_{good}(j)) \geq Dn^{-\frac{\alpha}{2}}\mathbb{P}(E_{good}(j)) \geq Dn^{-\frac{\alpha}{2}}c_0.\]
Thus from~(\ref{var_wl_one}), we have\[var(W_{n+1}) \geq \sum_{j=2}^{n+1} \mathbb{E}G^2_j \geq \sum_{j=2}^{n+1} (\mathbb{E}|G_j|)^2 \geq D^2c_0^2n^{1-\alpha},\] proving~(\ref{var_p_app}).

\begin{figure}[tbp]
\centering
\includegraphics[width=3in, trim= 20 300 50 110, clip=true]{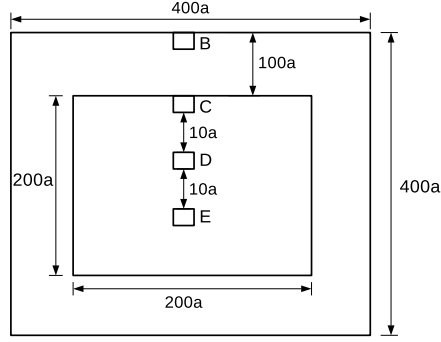}
\caption{The~\(a \times a\) squares~\(R^{B}_i,R^{C}_i,R^{D}_i\) and~\(R^{E}_i\) contained within the~\(400a \times 400a\) square~\(R_i\) with~\(a=\frac{A}{\sqrt{n}}.\) The~\(200a \times 200a\) square denotes~\(R^{small}_i\) and the centres of all the~\(a \times a\) squares
are equidistant from the left and right sides of~\(R_i.\)}
\label{good_square}
\end{figure}

\emph{Proof of~\((p1)-(p2)\)}: We recall the tiling of the unit square~\(S\) into small~\(\frac{400A}{\sqrt{n}} \times \frac{400A}{\sqrt{n}}\) squares~\(\{R_i\}.\) The side length of the smaller square~\(S_0 \subseteq S\) is~\(s_0,\) a constant and there are~\(W \geq \frac{1}{2}\cdot\left(\frac{s_0\sqrt{n}}{400A}\right)^2\) squares of~\(\{R_i\}\) completely contained in~\(S_0,\) which we label as~\(R_1,\ldots,R_W.\)  For each square~\(R_i,\) we let~\(R_i^{B},R_i^{C},R_i^{D}\) and~\(R_i^{E}\) be~\(\frac{A}{\sqrt{n}} \times \frac{A}{\sqrt{n}}\) subsquares contained in~\(R_i,\) as illustrated in Figure~\ref{good_square}, equidistant from the left and right sides of~\(R_i.\) Say that~\(R_i\) is a \(j-\)good square if the~\(\frac{A}{\sqrt{n}} \times \frac{A}{\sqrt{n}}\) square~\(R_i^{B}\) contains a single node of~\(\{X_k\}_{2 \leq k \neq j \leq n+1}\) and the rest of~\(R_i\) does not contain any node of~\(\{X_k\}_{2 \leq k \neq j \leq n+1}.\) Recalling that~\(R(1) \in \{R_i\}\) is the random square containing the node~\(X_1,\) we define the event
\begin{eqnarray}
E_{good}(j) &:=& \{R(1) \in \{R_i\}_{1 \leq i \leq W} \} \cap \{R(1) \text{ is } j-\text{good}\} \nonumber\\
&&\;\;\;\;\cap \{X_1 \in R^{C}(1)\} \cap \{X'_j \in R^{D}(1)\} \cap \{X_j \in R^E(1)\}.  \label{e_good_def}
\end{eqnarray}

Assuming the event~\(E_{good}(j)\) occurs, we now determine the change in the weight of the approximate overall TSP path  upon ``replacing" the node~\(X_j\) with the node~\(X'_j.\) The square~\(R(1)\) containing the node~\(X_1\) is contained within the square~\(S_0\) and so the weight of any edge with an endvertex in~\(R(1)\) is simply the Euclidean length of the edge (see statement of the Theorem). Since there are only two nodes~\(\{X_j,X_1\}\) present in the square~\(R^{small}(1) \subset R(1),\) the in-spanning path~\({\cal P}_{in}(n~+~1)\) is the single edge~\((X_j,X_1).\)

Next, the node~\(v\) present in the subsquare~\(R^B(1) \subset  R(1)\) is  \emph{the} node closest to but outside~\(R^{small}(1)\) and~\(v\) is closer to~\(X_1 \in R^C(1)\) than~\(X_j \in R^D(1).\) Therefore~\(v_{close} = v\) and the minimum length cross edge~\(e_{min} = (X_1,v).\) Combining we get
\begin{equation}\label{p_app_con}
{\cal P}_{app}(n+1) = {\cal P}_{app}\left(n+1,\{X_i\}_{1 \leq i \leq n+1}\right) = {\cal P}_{out}(n+1) \cup \{(X_1,X_j)\} \cup \{(X_1,v)\}.
\end{equation}
By an analogous analysis we have
\begin{equation}\label{p_app_con2}
{\cal P}_{app}\left(n+1,\{X_i\}_{1 \leq i \neq j\leq n+1} \cup X'_j\right) = {\cal P}_{out}(n+1) \cup \{(X_1,X'_j)\} \cup \{(X_1,v)\}.
\end{equation}

Denoting~\(W_{n+1}(X_j)\) and~\(W_{n+1}(X'_j)\) to be the respective weights of the approximate paths formed by~\(\{X_i\}_{1 \leq i \leq n+1}\)
and~\(\{X_i\}_{1 \leq i \neq j \leq n+1} \cup \{X'_j\},\) we get from~(\ref{p_app_con}) and~(\ref{p_app_con2}) that
\begin{eqnarray}
W_{n+1}(X_j)\ind(E_{good}(j)) &=& W_{n+1}(X'_j)\ind(E_{good}(j)) \nonumber\\
&&\;\;\;\;+\;\;\left(d^{\alpha}(X_1,X_j) - d^{\alpha}(X_1,X'_j)\right)\ind(E_{good}(j)) \nonumber\\
\label{one_diff_wl}
\end{eqnarray}
since~\(X_1,X_j\) and~\(X'_j\) belong to~\(S_0\) and so as metioned before, the weight of the edges is simply the corresponding Euclidean lengths raised
to the power~\(\alpha.\)

From Figure~\ref{good_square}, we have that the squares~\(R^{C}(1), R^{D}(1)\) and~\(R^{E}(1)\) are spaced~\(10a\) apart.
Since the nodes~\(X_1 \in R^{C}(1), X'_j \in R^{D}(1)\) and~\(X_j \in R^{E}(1),\) we have that~\(d(X_1,X_j)\geq d(X_1,X'_j) + 10a\)
and~\(10a \leq d(X_1,X'_j) \leq 14a.\) Thus~\[d^{\alpha}(X_1,X_j) = (d(X_1,X'_j) + 10a)^{\alpha} \geq (20a)^{\alpha}\]
and so we get from~(\ref{one_diff_wl}) that
\begin{equation}
W_{n+1}(X_j)\ind(E_{good}(j)) \geq W_{n+1}(X'_j)\ind(E_{good}(j))  + ((20a)^{\alpha}-(14a)^{\alpha})\ind(E_{good}(j)).
\label{one_diff_wl2}
\end{equation}

Setting~\(a = \frac{A}{\sqrt{n}}\) and using~\(A \geq 1,\) we finally get
\begin{eqnarray}
W_{n+1}(X_j)\ind(E_{good}(j)) &\geq& W_{n+1}(X'_j)\ind(E_{good}(j))  \nonumber\\
&&\;\;\;\;+\;\;((20)^{\alpha}-(14)^{\alpha})\cdot\left(\frac{1}{\sqrt{n}}\right)^{\alpha}\ind(E_{good}(j)),\nonumber\\
\label{one_diff_wl3}
\end{eqnarray}

It remains to estimate the probability of the event~\(E_{good}(j)\) and we use Poissonization. None of the constants appearing below depend on the index~\(j.\) As before, we let~\({\cal P}\) be a Poisson process with intensity~\(nf(.)\) in the unit square with corresponding probability measure~\(\mathbb{P}_0.\) A~\(\frac{400A}{\sqrt{n}} \times \frac{400A}{\sqrt{n}}\) square in~\(R_i\in \{R_k\}\) is good if the~\(\frac{A}{\sqrt{n}} \times \frac{A}{\sqrt{n}}\) subsquare~\(R^{B}_i\) contains exactly one node of~\({\cal P}\) and the rest of~\(R_i\) contains no node of~\({\cal P}.\) The mean number of nodes of~\({\cal P}\) in~\(R_i\) is~\(n\int_{R_i} f(x)dx\)
and~\(\mathbb{P}_0\left(R_i \text{ contains exactly one node of }{\cal P}\right)\) equals
\begin{equation}\label{p_r}
n\int_{R_i} f(x)dx \exp\left(-n\int_{R_i} f(x)dx\right) \geq A^2\epsilon_1e^{-A^2\epsilon_2} \geq \epsilon_1e^{-4\epsilon_2}
\end{equation}
where the first inequality in~(\ref{p_r}) follows from the bounds for~\(f(.)\) in~(\ref{f_eq}) and the second inequality in~(\ref{p_r}) is true since~\(1 \leq A \leq 1 + \frac{1}{\log{n}} \leq 2.\) Given that~\(R_i\) contains a single node of~\({\cal P},\) the node is distributed in~\(R_i\)
with density~\(f(.)\) and so is present in the subsquare~\(R^{B}_i\) with conditional probability~\(\frac{\int_{R^{B}_i}f(x)dx}{\int_{R_i}f(x)dx} \geq \frac{\epsilon_1}{(400)^2\epsilon_2},\) again using the bounds for~\(f(.)\) in~(\ref{f_eq}); i.e.,
\begin{equation}\label{cond_prob}
\mathbb{P}_0\left(R_i \text{ is good}|R_i \text{ contains exactly one node of }{\cal P}\right) = \frac{\int_{R^{B}_i}f(x)dx}{\int_{R_i}f(x)dx} \geq \frac{\epsilon_1}{(400)^2\epsilon_2}.
\end{equation}

Combining~(\ref{cond_prob}) and~(\ref{p_r}) we get that~\(\mathbb{P}_0(R_i \text{ is good}) \geq p_0\) for some constant~\(p_0 > 0\) not depending on~\(i.\) The Poisson process is independent on disjoint sets and so~\(N_{good} = \sum_{i=1}^{W} \ind(R_i \text{ is good})\) is a sum of independent Bernoulli random variables with mean at least~\(p_0,\) where we recall that~\(\{R_i\}_{1 \leq i \leq W}\) are all the squares contained completely inside the smaller square~\(S_0.\) Using the standard deviation estimate~(\ref{std_dev_down}) in Appendix with~\(m = W, \mu_1 =p_0\) and~\(\epsilon = \frac{1}{2},\)
we therefore have that~\(\mathbb{P}_0\left(N_{good} \geq \frac{p_0W}{2}\right) \geq 1-e^{-C_1W} \geq 1-e^{-2C_2n}\) for some positive constants~\(C_1\) and~\(C_2,\) since~\(W \geq \frac{1}{2}\cdot\left(\frac{s_0\sqrt{n}}{400A}\right)^2\) (see discussion prior to~(\ref{e_good_def})).

Using the dePoissonization formula~(\ref{de_poiss_ax}), we therefore get that
\[\mathbb{P}\left(N_{good} \geq \frac{p_0W}{2}\right) \geq 1-C_3\sqrt{n} e^{-2C_2n} \geq 1-e^{-C_2n}\] where~\(N_{good}\) is the number of~\(j-\)good squares (containing exactly one node of~\(\{X_k\}_{2 \leq k \neq j \leq n+1}\)) contained within~\(S_0.\) Given that~\(N_{good} \geq \frac{p_0W}{2} \geq C_4n,\) where~\(C_4 >0\) is a constant (see previous paragraph), the nodes~\(X_1,X'_j\) and~\(X'_j\) are all present in the same~\(\frac{400A}{\sqrt{n}} \times \frac{400A}{\sqrt{n}}\) good square~\(R_i \subset S_0\) in the corresponding~\(\frac{A}{\sqrt{n}} \times \frac{A}{\sqrt{n}}\) subsquares~\(R^{C}_i,R^{D}_i\) and~\(R^{E}_i\) with probability
at least~\(\left(\frac{p_0W}{2} \epsilon_1 \frac{A^2}{n}\right)^3 \geq p_1\) for some constant~\(p_1 >0,\) using the bounds for~\(f(.)\) in~(\ref{f_eq}).
From~(\ref{e_good_def}) we therefore get that~\(\mathbb{P}(E_{good}(j)) \geq p_1(1-e^{-C_2 n}).\)~\(\qed\)



\setcounter{equation}{0}
\renewcommand\theequation{\thesection.\arabic{equation}}
\section{Proof of convergence results}\label{pf_tsp_conv}
In this section, we obtain convergence results for~\(TSP_n\) appropriately scaled and centred.
We have the following result.
\begin{Theorem}\label{tsp_metric}
\((i)\) If~\(0 < \alpha < 2(\sqrt{2}-1),\) then~\[\frac{1}{n^{1-\frac{\alpha}{2}}}\left(TSP_n - \mathbb{E}TSP_n\right) \longrightarrow 0 \text{ a.s.} \] as~\(n \rightarrow \infty.\)\\
\((ii)\)Suppose the edge weight function~\(h\) is a metric and the edge weight exponent is~\(0 < \alpha  <1.\)
We have that~\[\frac{1}{n^{1-\frac{\alpha}{2}}}\left(TSP_n - \mathbb{E}TSP_n\right) \longrightarrow 0 \text{a.s.\ }\]
as~\(n \rightarrow \infty.\)
\end{Theorem}
The advantage of allowing~\(h\) to be a metric is that the TSP length satisfies the monotonocity property if the edge weight exponent is~\(\alpha \leq 1;\) in other words, adding more nodes increases the length of the TSP. We then obtain a.s.\ convergence via a subsequence argument by simply considering an upper bound on the TSP formed by the nodes with indices present outside the subsequence (see Section~\ref{pf_tsp_conv}).

To construct an example of~\(h\) which is a metric, we recall that~\(S = \left[-\frac{1}{2},\frac{1}{2}\right]^2\) is the unit square and define~\(h_1 : S \rightarrow (0,\infty)\)
as~\(h_1(u,v) = |(1+u)^2 -(1+v)^2|\) for~\(u,v \in \left[-\frac{1}{2},\frac{1}{2}\right]\) so that~\(h_1\) is a metric. Also,
\begin{equation}\label{bd1}
|u-v| \leq h_1(u,v) = |u-v|(2+u+v) \leq 3|u-v|.
\end{equation}
For~\(x=(x_1,x_2), y = (y_1,y_2),\) define
\[h(x,y) = h_1(x_1,y_1) + h_1(x_2,y_2)\] so that~\(h\) is a metric in~\(S\) and  using~(\ref{bd1}) we also get
\begin{equation}\nonumber
d(x,y) \leq |x_1-y_1| + |x_2-y_2| \leq h(x,y) \leq 3(|x_1-y_1| + |x_2-y_2|) \leq 3\sqrt{2}d(x,y).
\end{equation}

To prove Theorem~\ref{tsp_metric}, we have the following preliminary lemma.
We recall that~\(TSP_k\) denotes the length of the TSP cycle formed by the nodes~\(\{X_i\}_{1 \leq i \leq k}\)
as defined in~(\ref{min_weight_cycle}). For a constant~\(M > 0\) we let
\begin{equation}\label{jn_def}
J_n =  J_n(M) := \max_{n^M \leq k < (n+1)^M} \left|TSP_k - TSP_{n^M}\right|
\end{equation}
and have the following result.
\begin{Lemma}\label{j_lem} There is a constant~\(D > 0\) such that
\begin{equation}\label{jn_est}
\mathbb{E}\left(\frac{J_n}{n^{M\left(1-\frac{\alpha}{2}\right)}}\right) \leq  \frac{D}{n^{1-\frac{\alpha}{2}}} \text{ and } \mathbb{E}\left(\frac{J_n}{n^{M\left(1-\frac{\alpha}{2}\right)}}\right) ^2 \leq \frac{D}{n^{2-y}},
\end{equation}
where~\(y = \max\left(M\frac{\alpha}{2},\alpha\right).\) If the edge weight function~\(h\) is a metric and~\(0 < \alpha \leq 1,\) then
\begin{equation}\label{jn_est22}
\mathbb{E}\left(\frac{J_n}{n^{M\left(1-\frac{\alpha}{2}\right)}}\right) \leq  \frac{D}{n^{1-\frac{\alpha}{2}}} \text{ and } \mathbb{E}\left(\frac{J_n}{n^{M\left(1-\frac{\alpha}{2}\right)}}\right) ^2 \leq \frac{D}{n^{2-\alpha}}.
\end{equation}
\end{Lemma}

\emph{Proof of~(\ref{jn_est}) in Lemma~\ref{j_lem}}:  From the partial one node estimate~(\ref{tsp_up22}), we have that~\(TSP_{l} \leq TSP_{l+1} + g_{1,l},\)
where~\(g_{1,l} \geq 0\) is such that
\begin{equation}\label{ef_est2}
\max\left(\mathbb{E}g_{1,l}, \mathbb{E}g^2_{1,l}\right) \leq \frac{D}{l^{\frac{\alpha}{2}}}
\end{equation}
for some constant~\(D_1 > 0,\) using estimates~(\ref{g0_prop1}) and~(\ref{g1_prop1}).
For~\(n^{M}+1 \leq k < (n+1)^{M},\)
we therefore have~\[TSP_{n^{M}} \leq TSP_{k} + \sum_{l=n^{M}}^{k}g_{1,l}.\] To get an estimate in the reverse direction, we use relation~(\ref{mono_2}) in Appendix and get that~\[TSP_{k} \leq TSP_{n^{M}} + TSP(X_{n^{M}+1},\ldots,X_k)+ (c_2\sqrt{2})^{\alpha}.\] Combining we get~\(J_n \leq J_{n,1} + J_{n,2} + (c_2\sqrt{2})^{\alpha}\)
where
\[J_{n,1} := \sum_{l=n^{M}}^{(n+1)^M}g_{1,l} \text{ and } J_{n,2} := \max_{n^{M}+1 \leq k < (n+1)^{M}} TSP(X_{n^{M}+1},\ldots,X_k).\]

We evaluate~\(J_{n,1}\) and~\(J_{n,2}\) separately. Using~(\ref{ef_est2}), we have that
\begin{equation}\label{ejn1}
\mathbb{E}J_{n,1} \leq \sum_{l=n^{M}}^{(n+1)^M} \frac{D}{l^{\frac{\alpha}{2}}} \leq ((n+1)^{M}-n^{M})\frac{D}{n^{\frac{M\alpha}{2}}} \leq \frac{2MD}{n} \cdot n^{M\left(1-\frac{\alpha}{2}\right)}
\end{equation}
since
\begin{equation}\label{nm_est}
(n+1)^{M} -n^{M} \leq M(n+1)^{M-1} \leq 2Mn^{M-1}
\end{equation}
for all~\(n\) large. Similarly, using~\((\sum_{i=1}^{t}a_i)^2 \leq t\sum_{i=1}^{t} a_i^2,\) we get \[\mathbb{E}J^2_{n,1} \leq ((n+1)^{M}-n^{M})\sum_{l=n^{M}}^{(n+1)^M} \mathbb{E}g_{1,l}^2 \leq ((n+1)^{M}-n^{M})\sum_{l=n^{M}}^{(n+1)^M} \frac{D}{l^{\frac{\alpha}{2}}},\]
using~(\ref{ef_est2}). Again using~(\ref{nm_est}), we then get that~
\begin{equation}\label{ejn11}
\mathbb{E}J_{n,1}^2 \leq D_1 n^{2M-2-M\frac{\alpha}{2}}
\end{equation}
for some constant~\(D_1 > 0.\)

To estimate~\(J_{n,2},\) we let~\(L_k := TSP(X_{n^{M}+1},\ldots,X_k)\)
so that\\\(L_k \leq (k-n^{M}-1)\cdot(c_2 \sqrt{2})^{\alpha}\) by~(\ref{eq_met1}) since there are~\(k-n^{M}-1\) edges each of length at most~\(\sqrt{2}.\)
Thus
\begin{equation}\label{jn_crude}
J_{n,2} \leq ((n+1)^{M}-n^{M})\cdot(c_2 \sqrt{2})^{\alpha} \leq D n^{M-1}
\end{equation}
for some constant~\(D > 0\) using~(\ref{nm_est})  and if~\(x = (M-1)\left(1-\frac{\alpha}{2}\right),\)  then
\begin{equation}\label{tspk2}
\max_{n^{M}+1 \leq k < n^{M}+1+n^{x}} L_k \leq D_1n^{x} = D_1n^{(M-1)\left(1-\frac{\alpha}{2}\right)},
\end{equation}
for some constant~\(D_1 >0.\)  To bound~\(L_k\) for the remaining range of~\(k,\) we use the deviation estimate in Theorem~\ref{tsp_thm}. For~\(\Delta > 0,\) and~\(n^{M}+1+n^{x} \leq k \leq (n+1)^{M}\) we have from~(\ref{tsp_up_bounds}) that~\(\mathbb{P}\left(L_{k} \leq D_1(k-n^{M}-1)^{1-\frac{\alpha}{2}} \right) \geq 1-\frac{1}{k^{\Delta}}\) for some constant~\(D_1 > 0\) and so for~\(n^{M}+1+n^{x} \leq k < (n+1)^{M},\) we have that
\begin{equation}\label{tspk1}
L_{k} \leq D_1((n+1)^{M}-n^{M})^{\left(1-\frac{\alpha}{2}\right)} \leq D_2n^{(M-1)\left(1-\frac{\alpha}{2}\right)}
\end{equation}
for some constant~\(D_2 >0,\) with probability at least~\(1-\frac{1}{k^{\Delta}} \geq 1-\frac{1}{n^{M\Delta}},\) where the final estimate in~(\ref{tspk1}) follows from~(\ref{nm_est}). Setting~\(D_3 = \max(D_1,D_2),\) where~\(D_1\) and~\(D_2\) are as in~(\ref{tspk2}) and~(\ref{tspk1}), respectively, we therefore get that
\begin{equation}\label{jn2_prob}
\mathbb{P}\left(J_{n,2} \leq D_3 n^{(M-1)\left(1-\frac{\alpha}{2}\right)}\right) \geq 1-\frac{1}{n^{M\Delta}}((n+1)^{M}-n^{M}) \geq 1-\frac{2Mn^{M-1}}{n^{M\Delta}},
\end{equation}
by~(\ref{nm_est}).

From~(\ref{jn2_prob}) and~(\ref{jn_crude}), we get
\begin{equation}\label{ejn2}
\mathbb{E}J_{n,2} \leq D n^{(M-1)\left(1-\frac{\alpha}{2}\right)} + Dn^{M-1}\cdot\frac{Dn^{M-1}}{n^{M\Delta}} \leq 2Dn^{(M-1)\left(1-\frac{\alpha}{2}\right)},
\end{equation}
provided~\(\Delta > 0\) large. Choosing~\(\Delta > 0\) larger if necessary and arguing as above, we also get that
\begin{equation}\label{ejn22}
\mathbb{E}J_{n,2}^2 \leq 4D^2 n^{(M-1)(2-\alpha)}.
\end{equation}
From~(\ref{ejn1}) and~(\ref{ejn2}), we get the desired bound for~\(\mathbb{E}J_n\) in~(\ref{jn_est}).
Also, using~\((a+b)^2 \leq 2(a^2+b^2)\) and~(\ref{ejn11}) and~(\ref{ejn22}),
we get the desired bound for~\(\mathbb{E}J_n^2\) in~(\ref{jn_est}).~\(\qed\)

\emph{Proof of~(\ref{jn_est22}) in Lemma~\ref{j_lem}}: From~(\ref{mono_2}) and the monotonicity relation~(\ref{mono}) we get that
that~\(J_n \leq TSP(X_{n^{M}+1},\ldots,X_{(n+1)^{M}})\)
and so
\begin{equation}\label{jn_est222}
\mathbb{E}J_n^2 \leq \mathbb{E}TSP^2(X_{n^{M}+1},\ldots,X_{(n+1)^{M}}) \leq D_1\left((n+1)^{M}-n^{M}\right)^{2-\alpha} \leq D_2n^{(M-1)(2-\alpha)}
\end{equation}
for some constants~\(D_1,D_2 > 0,\) where the second inequality in~(\ref{jn_est222}) follows from the expectation upper bound for~\(TSP^{k}(.)\) in~(\ref{exp_tsp_bound}) with~\(k=2\) and the final inequality in~(\ref{jn_est222}) is true by~(\ref{nm_est}). Thus
\[\left(\frac{\mathbb{E}J_n}{n^{M\left(1-\frac{\alpha}{2}\right)}}\right)^2 \leq \frac{\mathbb{E}J^2_n}{n^{2M-M\alpha}} \leq \frac{D_2}{n^{M(2-\alpha)}},\]
proving~(\ref{jn_est22}).~\(\qed\)

\emph{Proof of Theorem~\ref{tsp_metric}~\((i)\)}: We prove a.s.\ convergence via a subsequence argument
using the variance upper bound for~\(TSP_n\) obtained in Theorem~\ref{tsp_var}. Letting~\(\epsilon_0 := \frac{(2-\alpha)(1+\alpha)}{2+\alpha},\) we use~\(\alpha < 2(\sqrt{2}-1) < 1\) to get that~\(\frac{2}{\alpha} > \frac{2+\alpha}{2-\alpha} = \frac{1}{2-\alpha-\epsilon_0}.\) We therefore choose~\(M>0\) such that~\(\frac{2}{\alpha} > M > \frac{1}{2-\alpha-\epsilon_0}\) and let~\(\epsilon > \epsilon_0\) be sufficiently close to~\(\epsilon_0\)
so that~\(\frac{2}{\alpha} > M > \frac{1}{2-\alpha-\epsilon}.\) Fixing such~\(M\) and~\(\epsilon\) we also have that~\(y = \max\left(M\frac{\alpha}{2},\alpha\right) < 1.\) For simplicity, we treat~\(n^{M}\) as an integer for all~\(n\) and work the subsequence~\(\{n^{M}\}.\)

From the variance estimate in Theorem~\ref{tsp_var}~\((i)\) we have that~\(var\left(\frac{TSP_n}{n^{1-\frac{\alpha}{2}}}\right) \leq \frac{D}{n^{2-\alpha-\epsilon}}\) for some constant~\(D >0\) and all~\(n\) large. Since~\(M(2-\alpha-\epsilon) > 1,\) we get by an application of the Borel-Cantelli Lemma that~\[\frac{1}{n^{M\left(1-\frac{\alpha}{2}\right)}}(TSP_{n^{M}} - \mathbb{E}TSP_{n^M}) \longrightarrow 0 \;\;\;\;a.s.\ \] as~\(n \rightarrow \infty.\) To prove convergence along the subsequence~\(a_n = n,\) we let~\(J_n\) be as in~(\ref{jn_def})
and obtain from~(\ref{jn_est}) that~\[\left(\frac{\mathbb{E}J_n}{n^{M\left(1-\frac{\alpha}{2}\right)}}\right)^2 \leq \frac{\mathbb{E}J^2_n}{n^{M(2-\alpha)}}   \leq \frac{D_1}{n^{2-y}} \longrightarrow 0\] as~\(n \rightarrow \infty\) and since~\(y < 1\) (see the first paragraph of this proof), we also get from Borel-Cantelli Lemma that~\(\frac{J_n}{n^{M(1-\frac{\alpha}{2})}} \longrightarrow 0\) a.s.\
as~\(n \rightarrow \infty.\) For~\(n^M \leq k < (n+1)^M\) we then write
\begin{eqnarray}
\frac{|TSP_k - \mathbb{E}TSP_{k}|}{k^{1-\frac{\alpha}{2}}} &\leq& \frac{|TSP_k - TSP_{n^M}|}{k^{1-\frac{\alpha}{2}}} + \frac{\mathbb{E}|TSP_{k} - TSP_{n^M}|}{k^{1-\frac{\alpha}{2}}} \nonumber\\
&\leq& \frac{J_n}{k^{1-\frac{\alpha}{2}}} + \frac{\mathbb{E}J_n}{k^{1-\frac{\alpha}{2}}} \nonumber\\
&\leq& \frac{J_n}{n^{M\left(1-\frac{\alpha}{2}\right)}} + \frac{\mathbb{E}J_n}{n^{M\left(1-\frac{\alpha}{2}\right)}} \label{tsp_as}
\end{eqnarray}
and get that~\(\frac{TSP_{k} - \mathbb{E}TSP_k}{k^{1-\frac{\alpha}{2}}} \longrightarrow 0\) a.s.\ as~\(n \rightarrow \infty.\)~\(\qed\)

\emph{Proof of Theorem~\ref{tsp_metric}~\((ii)\)}: Here we prove a.s.\ convergence for~\(0 < \alpha  \leq 1\) assuming that the edge weight function~\(h\) is a metric. The proof is analogous as before with some modifications. Let~\(2-\alpha > \epsilon >  (2-\alpha)\cdot \frac{1+\alpha}{2+\alpha}.\)
From the variance estimate in Theorem~\ref{tsp_var}~\((i)\) we have that~\[var\left(\frac{TSP_n}{n^{1-\frac{\alpha}{2}}}\right) \leq \frac{D}{n^{2-\alpha-\epsilon}}\] for some constant~\(D >0\) and all~\(n\) large. Letting~\(M > 0\) be a large integer so that~\(M(2-\alpha-\epsilon) > 1,\) we get by an application of the Borel-Cantelli Lemma that \[\frac{1}{n^{M\left(1-\frac{\alpha}{2}\right)}}(TSP_{n^{M}} - \mathbb{E}TSP_{n^M}) \longrightarrow 0 \text{ a.s.\ }\] as~\(n \rightarrow \infty.\)

To prove a.s.\ convergence along the sequence~\(a_n = n\) we let~\(J_n\) be as in~(\ref{jn_def}) and
get from~(\ref{jn_est22}) in property~\((p2)\) above that
\[\left(\frac{\mathbb{E}J_n}{n^{M\left(1-\frac{\alpha}{2}\right)}}\right)^2 \leq \frac{\mathbb{E}J^2_n}{n^{2M-M\alpha}} \leq \frac{D}{n^{M(2-\alpha)}} \longrightarrow 0\] as~\(n \rightarrow \infty,\) for some constant~\(D > 0.\) Thus~\(\frac{\mathbb{E}J_n}{n^{M\left(1-\frac{\alpha}{2}\right)}} \longrightarrow 0\) as~\(n \rightarrow \infty\) and since~\(M(2-\alpha) > M(2-\alpha-\epsilon) > 1,\) we get from Borel-Cantelli lemma that~\(\frac{J_n}{n^{M\left(1-\frac{\alpha}{2}\right)}} \longrightarrow 0\) a.s.\ as~\(n \rightarrow \infty.\) Arguing as in~(\ref{tsp_as}) above, we then get the desired a.s.\ convergence.~\(\qed\)

\setcounter{equation}{0}
\renewcommand\theequation{\thesection.\arabic{equation}}
\section{Uniform TSPs} \label{tsp_uni}
In this Section, we assume that the nodes~\(\{X_i\}_{1 \leq i \leq n}\) are uniformly distributed in the unit square and obtain bounds on the asymptotic values of the expected weight, appropriately scaled and centred. We assume that the positive edge weight function~\(h : \mathbb{R}^2 \times \mathbb{R}^2 \rightarrow \mathbb{R}\) satisfies~(\ref{eq_met1}) along with the following two properties:\\
\((b1)\) For every~\(a > 0\) we have
\[h(au,av) = a\cdot h(u,v) \text{ for all } u,v \in \mathbb{R}^2\]
\((b2)\) There exists a constant~\(h_0 > 0\) such that for all~\(b \in \mathbb{R}^2\) we have
\begin{equation}\label{h0_def}
h(b+u,b+v) \leq h_0 \cdot h(u,v).
\end{equation}
For example, recalling that~\(d(u,v)\) denotes the Euclidean distance between~\(u\) and~\(v,\) we have that the function
\begin{equation}\label{h_ex}
h(u,v) = d(u,v) + \frac{1}{2} |d(u,0) - d(v,0)|
\end{equation}
is a metric since
\[d(u,0) \leq d(0,v) + d(v,u) \text{ and } d(v,0) \leq d(0,u) + d(u,v)\] by triangle inequality and so
\[d(u,v) \leq h(u,v) \leq d(u,v) + \frac{1}{2} d(u,v).\] This implies that~\(h(u,v)\) satisfies~(\ref{eq_met1}) and by definition~\(h\) also
satisfies~\((b1).\) Moreover using the triangle inequality we have
\begin{eqnarray}
h(b+u,b+v) &=& d(b+u,b+v) + \frac{1}{2}|d(u+b,0)-d(v+b,0)|  \nonumber\\
&=& d(u,v) + \frac{1}{2}|d(u+b,0)-d(v+b,0)| \nonumber\\
&\leq& d(u,v) + \frac{1}{2} d(u+b,v+b) \nonumber\\
&=&d(u,v) + \frac{1}{2} d(u,v) \nonumber\\
&\leq& \frac{3}{2} h(u,v) \nonumber
\end{eqnarray}
by definition of~\(h\) in~(\ref{h_ex}). Thus~\((b2)\) is also satisfied with~\(h_0 = \frac{3}{2}.\)

We have the following result.
\begin{Theorem}\label{tsp_unif}
Suppose the distribution function~\(f(.)\) is uniform, i.e.,~\(\epsilon_1 = \epsilon_2 = 1\) in~(\ref{f_eq}) and the edge weight function~\(h(u,v)\) satisfies~(\ref{eq_met1}) and properties~\((b1)-(b2)\) above.  If the edge weight exponent is~\(0 < \alpha < 2,\) then
\begin{equation}\label{lim_sup_est}
0 < \liminf_{n} \frac{\mathbb{E} TSP_{n}}{n^{1-\frac{\alpha}{2}}} \leq \limsup_{n} \frac{\mathbb{E} TSP_{n}}{n^{1-\frac{\alpha}{2}}} \leq h_0^{\alpha} \cdot \liminf_{n} \frac{\mathbb{E} TSP_{n}}{n^{1-\frac{\alpha}{2}}} < \infty.
\end{equation}
\end{Theorem}
Thus the scaled weight of the TSP remains bounded within a factor of~\(h_0^{\alpha}.\)

Letting~\[J_n = \max_{n^{2} \leq k < (n+1)^{2}} |TSP_{k} - TSP_{n^{2}}|\] be as in~(\ref{jn_def}) with~\(M = 2,\) we get for~\(n^{2} \leq k < (n+1)^{2}\) that
\[\frac{\mathbb{E}TSP_k}{k^{1-\frac{\alpha}{2}}} \leq \frac{\mathbb{E}TSP_k}{n^{2-\alpha}} \leq \frac{\mathbb{E}TSP_{n^{2}}}{n^{2-\alpha}} + \frac{\mathbb{E}J_n}{n^{2-\alpha}}\]
and
\[\frac{\mathbb{E}TSP_k}{k^{1-\frac{\alpha}{2}}} \geq \frac{\mathbb{E}TSP_k }{(n+1)^{2-\alpha}} \geq \frac{\mathbb{E}TSP_{n^{2M}}}{(n+1)^{2-\alpha}}  - \frac{\mathbb{E} J_n}{(n+1)^{2-\alpha}}.\]
Using the estimate~(\ref{jn_est}) with~\(M = 2,\) we have~\(\frac{\mathbb{E}J_n}{n^{2-\alpha}} \longrightarrow 0\) as~\(n \rightarrow \infty\) and so we get
\[ \limsup_{n} \frac{\mathbb{E} TSP_{n}}{n^{1-\frac{\alpha}{2}}} = \limsup_n \frac{\mathbb{E} TSP_{n^2}}{n^{2-\alpha}}\]
and
\[ \liminf_{n} \frac{\mathbb{E} TSP_{n}}{n^{1-\frac{\alpha}{2}}} = \liminf_n \frac{\mathbb{E} TSP_{n^2}}{n^{2-\alpha}}.\]
In what follows therefore, we work with the sequence~\(a_n = n^2.\)

We begin with the following preliminary Lemma.\\
\begin{Lemma}\label{lem_thmu}
There is a constant~\(D > 0\) such that for any positive integers~\(n_1,n_2 \geq 1\) we have
\begin{eqnarray}
\mathbb{E} TSP_{n_1+n_2} \leq \mathbb{E}TSP_{n_1} + Dn_2^{1-\frac{\alpha}{2}} + 2(\sqrt{2})^{\alpha}. \label{tsp_ab2}
\end{eqnarray}
Also for any fixed integer~\(m \geq 1\) we have
\begin{equation}\label{km_su}
\limsup_n \frac{\mathbb{E} TSP_{n^{2}}}{n^{2-\alpha}} \leq \limsup_k \frac{\mathbb{E}TSP_{k^2m^2}}{(km)^{2-\alpha}}.
\end{equation}
\end{Lemma}

\emph{Proof of~(\ref{tsp_ab2}) in Lemma~\ref{lem_thmu}}: From the additive relation~(\ref{mono_2}) in Appendix we have
\[TSP_{n_1+n_2} \leq TSP_{n_1} + TSP(X_{n_1+1},\ldots,X_{n_1+n_2}) + 2(\sqrt{2})^{\alpha}\]
and using the TSP expectation upper bound in~(\ref{exp_tsp_bound}), the middle term
is bounded above by~\(Dn_2^{1-\frac{\alpha}{2}}\) for some constant~\(D > 0.\)~\(\qed\)


\emph{Proof of~(\ref{km_su}) in Lemma~\ref{lem_thmu}}: Fix an integer~\(m \geq 1\) and write~\(n = qm + s\)
where~\(q = q(n) \geq 1\) and~\(0 \leq s = s(n) \leq m-1\) are integers.
As~\(n \rightarrow \infty,\)
\begin{equation}\label{kn_def}
q(n) \longrightarrow \infty \text{ and }\frac{n}{q(n)} \longrightarrow m.
\end{equation}
Also using~\(s <m,\) we have~\(qm+s < (q+1)m\) and so~\(q_1 := (qm+s)^{2} - (qm)^{2} \) is bounded above by
\begin{equation}\label{q1_est}
q_1 \leq m^{2}\left((q+1)^{2}-q^{2}\right) \leq 4m^{2}q.
\end{equation}

Using property~\((b1)\) we therefore have that~\(\mathbb{E}TSP_{n^{2}} = \mathbb{E}TSP_{(qm+s)^{2}}\) equals
\[\mathbb{E}TSP_{(qm)^{2} + q_1} \leq \mathbb{E}TSP_{(qm)^{2}} + Dq_1^{1-\frac{\alpha}{2}} + 2(\sqrt{2})^{\alpha}\]
for some constant~\(D >0.\) From~(\ref{q1_est}), we have~\(Dq_1^{1-\frac{\alpha}{2}} \leq \frac{D_1}{q^{1-\frac{\alpha}{2}}}\cdot (qm)^{2-\alpha}\)
for some constant~\(D_1  >0\) and so
\begin{eqnarray}
\limsup_n \frac{\mathbb{E}TSP_{n^{2}}}{n^{2-\alpha}} &\leq& \limsup_n \left(\frac{qm}{n}\right)^{2-\alpha} \frac{\mathbb{E} TSP_{(qm)^{2}}}{(qm)^{2-\alpha}} \nonumber\\
&&\;\;\;\;\;\;\;\;\;+\;\;\limsup_n \left(\frac{qm}{n}\right)^{2-\alpha} \left(\frac{D_1}{q^{1-\frac{\alpha}{2}}} + \frac{2(\sqrt{2})^{\alpha}}{(qm)^{2-\alpha}}\right).\nonumber\\
\label{tsp_temp_km1}
\end{eqnarray}
Since~\(q(n) \longrightarrow \infty\) and~\(\frac{n}{q(n)} \longrightarrow m\) (see~(\ref{kn_def})),~\(m\) is fixed and~\(0 < \alpha < 2,\) the second term in~(\ref{tsp_temp_km1}) is zero and the first term in~(\ref{tsp_temp_km1}) equals~\(\limsup_n \frac{\mathbb{E} TSP_{(qm)^{2}}}{(qm)^{2-\alpha}}.\)

But since~\(q(n) \geq \frac{n-m}{m} \geq \frac{l-m}{m}\) for~\(n \geq l,\) we have that
\begin{equation}\label{gtu}
\sup_{n \geq l} \frac{\mathbb{E} TSP_{q^2m^2}}{(qm)^{2-\alpha}} = \sup_{n \geq l} \frac{\mathbb{E} TSP_{q^2(n)m^2}}{(q(n)m)^{2-\alpha}} \leq \sup_{k \geq \frac{l-m}{m}} \frac{\mathbb{E} TSP_{k^2m^2}}{(km)^{2-\alpha}}
\end{equation}
and as~\(l \uparrow \infty,\) the final term in~(\ref{gtu}) converges to the second term in~(\ref{km_su}).~\(\qed\)


If~\(\lambda := \liminf_n \frac{\mathbb{E}TSP_{n^{2}}}{n^{2-\alpha}}\) then from~(\ref{exp_tsp_bound}), we have that~\(\lambda > 0.\) We show below that for every~\(\epsilon > 0,\) there exists~\(m\) sufficiently large such that
\begin{equation}\label{second_step_est}
\limsup_{k} \frac{\mathbb{E} TSP_{k^2m^2}}{(km)^{2-\alpha}} \leq h_0^2 \cdot \lambda + \epsilon. 
\end{equation}
Since~\(\epsilon > 0\) is arbitrary, the desired bound in Theorem~\ref{tsp_unif} follows from~(\ref{km_su}). In the rest of this Section, we prove~(\ref{second_step_est}) via a series of Lemmas.

Let~\(k\) and~\(m\) be positive integers and distribute~\(k^2m^2\) nodes\\\(\{X_i\}_{1 \leq i \leq k^2m^2}\) independently
and uniformly in the unit square~\(S.\) Also divide~\(S\) into~\(k^{2}\) disjoint
squares~\(\{W_j\}_{1 \leq j \leq k^{2}}\) each of size~\(\frac{1}{k} \times \frac{1}{k}\) as in Figure~\ref{fig_squares}
so that the top left most square is labelled~\(W_1,\) the square below~\(W_1\) is~\(W_2\) and so on until
we reach the square~\(W_k\) intersecting the bottom edge of the unit square~\(S.\) The square to the right of~\(W_k\)
is then labelled~\(W_{k+1}\) and the square above~\(W_{k+1}\) is~\(W_{k+2}\) and so on.


The next result obtains an upper bound for~\(TSP_{k^2m^2}\) in terms of the local TSPs of the squares~\(W_j.\) For~\(1 \leq j \leq k^2,\) Let~\(N(j)\) be the number of nodes of~\(\{X_i\}\) present in the square~\(W_j\) and let~\(TSP_j(N(j))\) be the TSP length of the nodes present in~\(W_j.\) We have the following Lemma.
\begin{Lemma}\label{shift_tsp} There exists a constant~\(D > 0\) not depending on~\(k\) or~\(m\) such that
\begin{equation}\label{tsp_km_m3}
\mathbb{E}TSP_{k^2m^2} \leq \sum_{j=1}^{k^2} \mathbb{E}TSP_j(N(j)) + D\cdot k^{2-\alpha} + (c_2\sqrt{2})^{\alpha},
\end{equation}
where~\(c_2\) is as in~(\ref{eq_met1}).
\end{Lemma}
\emph{Proof of Lemma~\ref{shift_tsp}}: We first connect all nodes within each square~\(W_j\) by a path of minimum weight
to get local minimum weigh paths. We then join all these paths together  to get an overall spanning path containing all the~\(k^2m^2\) nodes. Joining the endvertices of this new path by an edge, we get a spanning cycle whose weight is at least~\(TSP_{k^2m^2}.\)

For~\(1 \leq j \leq k^{2},\) let~\({\cal P}(j)\) and~\({\cal C}(j)\) respectively denote the minimum spanning path
and minimum spanning cycle, containing all the~\(N(j)\) nodes. If~\(N(j) =0\) we set~\({\cal P}(j) = \emptyset\)
and if~\(N(j) \leq 2,\) we set~\({\cal P}(j) = {\cal C}(j).\) For any~\(1 \leq j \leq k^2m^2\)
we first show that
\begin{equation}\label{pj_cj_comp}
W({\cal P}(j)) \leq W({\cal C}(j)) \leq W({\cal P}(j)) + \left(\frac{c_2\sqrt{2}}{k}\right)^{\alpha}.
\end{equation}
To see that~(\ref{pj_cj_comp}) is true, let~\({\cal C}\) be any spanning cycle containing all the nodes
in the square~\(W_j.\) Removing an edge of~\({\cal C}\) we get a spanning path~\({\cal P}\)
and so~\(W({\cal P}(j)) \leq W({\cal P}) \leq W({\cal C}).\) Taking minimum over all spanning cycles
we obtain the first relation in~(\ref{pj_cj_comp}). Taking the minimum spanning path~\({\cal P}(j)\) with endvertices~\(u,v\)
and adding the edge~\((u,v),\) we obtain a spanning cycle~\({\cal C}_{new}\) and since the length of the edge~\((u,v)\)
is at most~\(\frac{\sqrt{2}}{k}\) we use~(\ref{eq_met1}) to get that
\[W({\cal C}(j)) \leq W({\cal C}_{new}) \leq W({\cal P}(j)) + \left(\frac{c_2\sqrt{2}}{k}\right)^{\alpha},\]
proving~(\ref{pj_cj_comp}).

For~\(1 \leq j \leq k^{2}-1,\) let~\(j~+~T^{next}_j := \min\{i \geq j+1 : W_i \text{ is not empty}\}\) be the next
nonempty square after~\(W_j,\) containing at least one node of~\(\{X_i\}_{1 \leq i \leq k^2m^2}.\) If no such square exists, set~\(j+T^{next}_j~=~k^2.\) If both~\(W_j\) and~\(W_{j+T^{next}_j}\) are nonempty, let~\(e_j\) be the edge with one endvertex in~\(W_j\) and the other endvertex in~\(W_{j+T^{next}_j},\)
having the smallest Euclidean length~\(d(e_j)\) and set~\(e_j = \emptyset\) and~\(d(e_j) = 0,\) otherwise.  The union~\[{\cal P}_{ov} := \cup_{1 \leq j \leq k^{2}} {\cal P}(j) \cup \cup_{1 \leq j \leq k^{2}-1} \{e_j\}\] is a spanning path containing all the~\(k^2m^2\) nodes and joining the endvertices of~\({\cal P}_{ov}\) (by an edge of length at most~\(\sqrt{2}\)), we again use~(\ref{eq_met1}) to get
\begin{eqnarray}
TSP_{k^2m^2} &\leq& W({\cal P}_{ov}) + (c_2\sqrt{2})^{\alpha} \nonumber\\
&=& \sum_{j=1}^{k^{2}} W({\cal P}(j))  + \sum_{j=1}^{k^{2}-1} h^{\alpha}(e_j)   +(c_2\sqrt{2})^{\alpha}\nonumber\\
&\leq& \sum_{j=1}^{k^{2}} W({\cal P}(j))  + c_2^{\alpha}\sum_{j=1}^{k^{2}-1} d^{\alpha}(e_j)   +(c_2\sqrt{2})^{\alpha}\nonumber\\
&\leq& \sum_{j=1}^{k^{2}} W({\cal C}(j))  + c_2^{\alpha}\sum_{j=1}^{k^{2}-1} d^{\alpha}(e_j) + (\sqrt{2})^{\alpha}, \label{tsp_km22}
\end{eqnarray}
using~(\ref{pj_cj_comp}).

Taking expectations in~(\ref{tsp_km22}) and recalling that~\(TSP_j(N(j)) := W({\cal C}(j))\) we get
\begin{equation}\label{tsp_km_2m}
\mathbb{E}TSP_{k^2m^2} \leq \sum_{j=1}^{k^2} \mathbb{E}TSP_j(N(j)) + c_2^{\alpha} \cdot k^{2} \max_{1 \leq j \leq k^{2}-1} \mathbb{E}d^{\alpha}(e_j) + (c_2\sqrt{2})^{\alpha}.
\end{equation}
We show that there exists a constant~\(D > 0\) not depending on~\(k\) or~\(m\) such that
\begin{equation}\label{s_eval}
\max_{1 \leq j \leq k^{2}-1} \mathbb{E}d^{\alpha}(e_j) \leq \frac{D}{k^{\alpha}}
\end{equation}
and this proves Lemma~\ref{shift_tsp}.
Indeed, the Euclidean length of the edge~\(e_j\) is at most~\(\frac{2T^{next}_j\sqrt{2}}{k}\)
and so~\(\mathbb{E}d^{\alpha}(e_j) \leq \frac{(2\sqrt{2})^{\alpha}}{k^{\alpha}} \mathbb{E}(T^{next}_j)^{\alpha}.\) Next, for any~\(l \geq 1\) we have~\(T^{next}_j > l\) if and only if~\(W_{j+1},\ldots,W_{j+l}\) are empty, which happens with probability~\(\left(1-l\cdot\frac{1}{k^2}\right)^{k^2m^2} \leq e^{-lm^{2}} \leq e^{-l},\) since~\(m \geq 1.\) Thus
\[\mathbb{E}(T^{next}_j)^{\alpha} \leq \mathbb{E}(T^{next}_j)^{2} \leq \sum_{l \geq 1} l \mathbb{P}\left(T^{next}_j \geq l\right) \leq \sum_{l \geq 1} l e^{-(l-1)} \leq \frac{1}{(1-e^{-1})^2},\] proving~(\ref{s_eval}).~\(\qed\)

Lemma~\ref{shift_tsp} is useful in the following manner. For~\(1 \leq j \leq k^2\) let~\(TSP(N(j))\) be the \emph{shifted} TSP of the nodes present in~\(W_j\) defined as follows: If~\(v_1,\ldots,v_t\) are the nodes of~\(\{X_i\}\) present in the square~\(W_j\) with centre~\(s_j,\) then~\(TSP(N(j))\) is the TSP formed by the nodes~\(v_1-s_j,\ldots,v_t-s_j.\) From the translation relation~(\ref{trans_rel}) in Appendix, we have~\(TSP_j(N(j)) \leq h_0^{\alpha} \cdot TSP(N(j))\) and so we get from Lemma~\ref{shift_tsp} that
\begin{eqnarray}\label{cruc_rec}
\mathbb{E}TSP_{k^2m^2} &\leq& h_0^{\alpha} \cdot \sum_{j=1}^{k^2} \mathbb{E}TSP(N(j)) + D\cdot k^{2-\alpha} + (c_2\sqrt{2})^{\alpha}\nonumber\\
&=& h_0^{\alpha} \cdot k^2 \cdot \mathbb{E}TSP(N(1)) + D \cdot k^{2-\alpha} + (c_2\sqrt{2})^{\alpha}\nonumber
\end{eqnarray}
To evaluate~\(TSP(N(1)),\) we let~\(0 < \gamma < \frac{1}{2}\) and~\(\frac{(2-\alpha)(1+\alpha)}{2+\alpha} < \epsilon_1 <2-\alpha\) be constants (not depending on~\(k\) or~\(m\)) such that
\begin{equation}\label{gam_choice}
2-\alpha + 2\gamma - 2\epsilon_1 > 0
\end{equation}
This is possible since~\( 2\frac{(2-\alpha)(1+\alpha)}{2+\alpha} -(2-\alpha) = \frac{\alpha(2-\alpha)}{2+\alpha} < 1\) for all~\(0 < \alpha < 2\)
and so we choose~\(\epsilon_1\) greater than but sufficiently close to~\(\frac{(2-\alpha)(1+\alpha)}{2+\alpha}\) and~\(2\gamma\) less than but sufficiently close to one so that~(\ref{gam_choice}) holds.


We have the following Lemma.
\begin{Lemma}\label{i1_i2_lem} There is a constant~\(D > 0\) not depending on~\(k\) or~\(m\) such that
\begin{equation}\label{tsp_n1}
\mathbb{E}TSP(N(1)) \leq  \frac{1}{k^{\alpha}} \left(\mathbb{E} TSP_{m^{2}} + D m^{2(\epsilon_1-\gamma)} + D \frac{m^{2-\alpha} }{m^{\alpha(1-2\gamma)}}\right)
\end{equation}
\end{Lemma}
\emph{Proof Lemma~\ref{i1_i2_lem}}: We write~\(\mathbb{E} TSP(N(1)) = I_1 + I_2,\) where~\[I_1 =\mathbb{E} TSP(N(1)) \ind(F_1), I_2 = \mathbb{E}TSP(N(1)) \ind(F_1^c)\] and~\(F_1 := \{m^{2}\left(1 - \frac{1}{m^{2\gamma}}\right)  \leq N(1) \leq m^{2} \left(1+ \frac{1}{m^{2\gamma}}\right)\}.\)
Each node~\(X_i, 1 \leq i \leq k^2m^2\) has a probability~\(\frac{1}{k^{2}}\) of being present in the~\(\frac{1}{k} \times \frac{1}{k}\) square~\(W_1.\) Therefore the number of nodes~\(N(1)\) in the square~\(W_1\) is binomially distributed with mean~\(\mathbb{E}N(1) = m^{2}\) and~\(var(N_1) \leq m^{2}.\) We therefore get from Chebychev's inequality that
\begin{equation}\label{fm_k_est}
\mathbb{P}(F^c_1)  \leq \frac{1}{m^{2(1-2\gamma)}}. 
\end{equation}

\emph{\underline{Evaluation of~\(I_1\)}}: We write~\[I_1 = \sum_{j = j_{low}}^{j_{up}} \mathbb{E} TSP(N(1)) \ind(N(1) = j),\]
where~\(j_{low} := m^{2}\left(1 - \frac{1}{m^{2\gamma}}\right)   \leq m^{2}\left(1 + \frac{1}{m^{2\gamma}}\right)   =: j_{up}.\) Given~\(N_1 = j,\) the nodes in~\(W_1\) are independently and uniformly distributed in~\(W_1\) and we let~\(\mathbb{E}TSP\left(j;\frac{1}{k}\right)\) be the
expected length of the TSP containing~\(j\) nodes independently and uniformly distributed in the~\(\frac{1}{k} \times \frac{1}{k}\) square centred at origin.
From the scaling relation~(\ref{tsp_scale}) in Appendix we have
\begin{equation}
I_1 = \sum_{j =j_{low}}^{j_{up}} \mathbb{E}TSP\left(j;\frac{1}{k}\right) \mathbb{P}(N(1) = j) = \frac{1}{k^{\alpha}} \sum_{j=j_{low}}^{j_{up}} \left(\mathbb{E}TSP_j \right)\mathbb{P}(N(1) = j), \label{temp1}
\end{equation}
by~(\ref{tsp_scale}).

Letting~\(2-\alpha > \epsilon_1 > \frac{(2-\alpha)(1+\alpha)}{2+\alpha}\) be as in~(\ref{gam_choice}), we have from the one node difference estimate~(\ref{cruc_tsp_est}) that for any integers~\(j_1\) and~\(j_2\) lying between~\(j_{low}\) and~\(j_{up},\)
the term~\(\mathbb{E}|TSP_{j_2} - TSP_{j_1}| \) is bounded above by
\begin{equation}
\sum_{u=j_{low}}^{j_{up}-1} \mathbb{E}|TSP_{u+1} - TSP_{u}| \leq \sum_{u=j_{low}}^{j_{up}-1} \frac{D_1}{u^{1-\epsilon_1}} \leq \frac{D_2(j_{up}-j_{low})}{m^{2(1-\epsilon_1)}} = 2D_2\frac{m^{2(1-\gamma)}}{m^{2(1-\epsilon_1)}},\label{temp3}
\end{equation}
where~\(D_1,D_2 > 0\) are constants not depending on~\(j_1\) or~\(j_2.\) Setting~\(j_1= m^{2}\) and~\(j_2 = j\) and using~(\ref{temp3}) we get~\(\mathbb{E}TSP_j \leq \mathbb{E}TSP_{m^{2}} + D_2 m^{2(\epsilon_1-\gamma)}\) for all~\(j_{low} \leq j \leq j_{up}.\) From the expression for~\(I_1\) in~(\ref{temp1}) we therefore have that
\begin{equation}
I_1 \leq \frac{1}{k^{\alpha}} \left(\mathbb{E} TSP_{m^{2}} + D_2 m^{2(\epsilon_1-\gamma)}\right). \label{i1_est_a}
\end{equation}

\underline{\emph{Evaluation of~\(I_2\)}}: There are~\(N(1)\) nodes in the~\(\frac{1}{k} \times \frac{1}{k}\) square~\(W_1\) and given~\(N(1) = l,\) the~\(l\) nodes are uniformly distributed in~\(W_1\) and so
\[\mathbb{E}(TSP(N(1))|N(1)=l) = \mathbb{E}TSP\left(l;\frac{1}{k}\right) = \frac{1}{k^{\alpha}}\mathbb{E}TSP_l \leq D_1\frac{l^{1-\frac{\alpha}{2}}}{k^{\alpha}}\]
for some constant~\(D_1 > 0\) using the expectation upper bound in~(\ref{exp_tsp_bound}).
Thus~
\begin{eqnarray}
\mathbb{E}TSP(N(1))\ind(N(1) = l) &=& \mathbb{E}(TSP(N(1))|N(1)=l)\mathbb{P}(N(1) = l) \nonumber\\
&\leq& D_1\frac{l^{1-\frac{\alpha}{2}}}{k^{\alpha}}\mathbb{P}(N(1)= l) \nonumber
\end{eqnarray}
and consequently~\(I_2 = \mathbb{E}TSP(N(1))\ind(F_1^c)\) equals
\begin{eqnarray}
&&\sum_{l\leq j_{low}} + \sum_{l \geq j_{up}} \mathbb{E}TSP(N(1))\ind(N(1) = l) \nonumber\\
&&\leq\;\;\; \sum_{l\leq j_{low}} + \sum_{l \geq j_{up}} D_1\frac{l^{1-\frac{\alpha}{2}}}{k^{\alpha}}\mathbb{P}(N(1)= l) \nonumber\\
&&=\;\;\; \frac{D_1}{k^{\alpha}} \mathbb{E}(N(1))^{1-\frac{\alpha}{2}}\ind(F_1^c). \label{i2_alt_est}
\end{eqnarray}

Using H{o}lder's inequality~\(\mathbb{E}XY \leq \left(\mathbb{E}X^{p}\right)^{\frac{1}{p}} \left(\mathbb{E}Y^q\right)^{\frac{1}{q}}\)
with~\(X = (N(1))^{1-\frac{\alpha}{2}},\)\\\(Y = \ind(F_1^c),p =\frac{2}{2-\alpha} >1\) and~\(q = \frac{2}{\alpha} > 1,\) we get
\begin{eqnarray}
\mathbb{E}(N(1))^{1-\frac{\alpha}{2}}\ind(F_1^c) &\leq& \left(\mathbb{E}N(1)\right)^{1-\frac{\alpha}{2}} \left(\mathbb{P}(F_1^c)\right)^{\frac{\alpha}{2}} \nonumber\\
&=& m^{2-\alpha}\left(\mathbb{P}(F_1^c)\right)^{\frac{\alpha}{2}} \nonumber\\
&\leq& m^{2-\alpha} \frac{1}{m^{\alpha(1-2\gamma)}},
\end{eqnarray}
by~(\ref{fm_k_est}). From~(\ref{i2_alt_est}) we therefore get that~
\begin{equation}\label{i2_est_b}
I_2 \leq \frac{D_2m^{2-\alpha}}{k^{\alpha}} \frac{1}{m^{\alpha(1-2\gamma)}}.
\end{equation}

From the estimate for~\(I_1\) in~(\ref{i1_est_a}), we therefore get that
\[\mathbb{E}TSP(N(1)) = I_1+I_2 \leq  \frac{1}{k^{\alpha}} \left(\mathbb{E} TSP_{m^{2}} + D_2 m^{2(\epsilon_1-\gamma)} + D_2 \frac{m^{2-\alpha} }{m^{\alpha(1-2\gamma)}}\right),\] proving Lemma~\ref{shift_tsp}.~\(\qed\)

Substituting~(\ref{tsp_n1}) of Lemma~\ref{shift_tsp} into~(\ref{tsp_km_m3}) we get
\begin{eqnarray}
\mathbb{E} TSP_{k^2m^2} &\leq& k^{2-\alpha}\left( h_0^{\alpha} \cdot \mathbb{E} TSP_{m^{2}} + D m^{2(\epsilon_1-\gamma)}  + D  \frac{m^{2-\alpha}}{m^{\alpha(1-2\gamma)}}+D\right) \nonumber\\
\label{tsp_km_m4}
\end{eqnarray}
for some constant~\(D > 0.\) Thus
\begin{equation}\nonumber
\limsup_k \frac{\mathbb{E} TSP_{k^2m^2}}{(km)^{2-\alpha}} \leq h_0^{\alpha }\cdot \frac{\mathbb{E} TSP_{m^{2}}}{m^{2-\alpha}} + \frac{D}{m^{2-\alpha+2\gamma-2\epsilon_1}} + \frac{D}{m^{\alpha(1-2\gamma)}} + \frac{D}{m^{2-\alpha}}
\end{equation}
for all~\(m\) large. By choice of~\(\gamma,\epsilon_1 > 0\) in~(\ref{gam_choice}) we have~\(2\gamma <1\)
and\\\(2-\alpha+2\gamma-2\epsilon_1 > 0\) and so
\begin{equation}
\limsup_k \frac{\mathbb{E} TSP_{k^2m^2}}{(km)^{2-\alpha}} \leq h_0^{\alpha} \cdot \frac{\mathbb{E} TSP_{m^{2}}}{m^{2-\alpha}} + \epsilon,\label{km_sw}
\end{equation}
for all~\(m\) large by definition of~\(\lambda\) in~(\ref{km_su}).  Letting~\(\{m_j\}\) be any sequence such that~\(\frac{\mathbb{E} TSP_{m^2_j}}{m_j^{2-\alpha}} \longrightarrow \lambda = \liminf_n \frac{\mathbb{E}TSP_n^2}{n^{2-\alpha}}\) and allowing~\(m \rightarrow \infty\) through the sequence~\(\{m_j\}\)  in~(\ref{km_sw}), we get that~\(\limsup_n \frac{\mathbb{E}TSP_n^2}{n^{2-\alpha}} \leq h_0^{\alpha} \cdot \lambda + \epsilon.\) Since~\(\epsilon  >0\) is arbitrary, this obtains~(\ref{second_step_est}).\(\qed\)

\setcounter{equation}{0}
\renewcommand\theequation{\thesection.\arabic{equation}}
\section{Appendix : Miscellaneous results}\label{appendix}
\subsection*{\em Standard deviation estimates}
We use the following standard deviation estimates for sums of independent Poisson and Bernoulli random variables.
\begin{Lemma}\label{app_lem}
Suppose~\(W_i, 1 \leq i \leq m\) are independent Bernoulli random variables satisfying~\(\mu_1 \leq \mathbb{P}(W_1=1) = 1-\mathbb{P}(W_1~=~0) \leq \mu_2.\) For any~\(0 < \epsilon < \frac{1}{2},\)
\begin{equation}\label{std_dev_up}
\mathbb{P}\left(\sum_{i=1}^{m} W_i > m\mu_2(1+\epsilon) \right) \leq \exp\left(-\frac{\epsilon^2}{4}m\mu_2\right)
\end{equation}
and
\begin{equation}\label{std_dev_down}
\mathbb{P}\left(\sum_{i=1}^{m} W_i < m\mu_1(1-\epsilon) \right) \leq \exp\left(-\frac{\epsilon^2}{4}m\mu_1\right)
\end{equation}
Estimates~(\ref{std_dev_up}) and~(\ref{std_dev_down}) also hold if~\(\{W_i\}\) are independent Poisson random variables with~\(\mu_1 \leq \mathbb{E}W_1 \leq \mu_2.\)
\end{Lemma}
For completeness, we give a quick proof.\\
\emph{Proof of Lemma~\ref{app_lem}}: First suppose that~\(\{W_i\}\) are independent Poisson with\\\(\mu_1 \leq \mathbb{E}W_i \leq \mu_2\) so that~\(\mathbb{E}e^{sW_i} = \exp\left(\mathbb{E}W_i (e^{s}-1)\right) \leq \exp\left(\mu_2(e^{s}-1)\right)\) for~\(s~>~0.\) By Chernoff bound
we then have \[\mathbb{P}\left(\sum_{i=1}^{m} W_i > m\mu_2(1+\epsilon)\right) \leq e^{-sm\mu_2(1+\epsilon)} \exp\left(m\mu_2(e^{s}-1)\right) = e^{m\mu_2 \Delta_1},\]
where~\(\Delta_1 = e^{s} -1 -s-s\epsilon.\) For~\(s \leq 1,\) we have the bound~\[e^{s}-1-s = \sum_{k\geq 2} \frac{s^{k}}{k!} \leq s^2\sum_{k\geq 2} \frac{1}{k!} = s^2(e-2) \leq s^2\] and so we set~\(s = \frac{\epsilon}{2}\) to get that~\(\Delta_1 \leq s^2-s\epsilon =\frac{-\epsilon^2}{4}.\)

Similarly for~\(s > 0,\) we have~\[\mathbb{E}e^{-sW_i} = \exp\left(\mathbb{E}W_i(e^{-s}-1)\right) \leq \exp\left(\mu_1(e^{-s}-1)\right)\] and so
\[\mathbb{P}\left(\sum_{i=1}^{m} W_i < m\mu_1(1-\epsilon)\right) \leq e^{sm\mu_1(1-\epsilon)} \exp\left(m\mu_1(e^{-s}-1)\right) = e^{-m\mu_1\Delta_2},\]
where~\(\Delta_2 = 1 -s-e^{-s} + s\epsilon.\) For~\(s \leq 1,\) the term~\(e^{-s}\leq 1-s + \frac{s^2}{2}\) and so we get~\(\Delta_2 \geq -\frac{s^2}{2}+s\epsilon =\frac{\epsilon^2}{2}\) for~\(s = \epsilon.\)

The proof for the Binomial distribution follows from the fact that if~\(\{W_i\}_{1 \leq i \leq m}\) are independent Bernoulli distributed with~\(\mu_1 \leq \mathbb{E}W_i \leq \mu_2,\) then for~\(s > 0\) we have~\(\mathbb{E}e^{s W_i} = 1-\mathbb{E}W_i + e^{s} \mathbb{E}W_i \leq \exp\left(\mathbb{E}W_i (e^{s}-1)\right) \leq \exp\left(\mu_2(e^{s}-1)\right)\) and similarly~\(\mathbb{E}e^{-s W_i}  \leq \exp\left(\mu_1(e^{-s}-1)\right).\) The rest of the proof is then as above.~\(\qed\)

\subsection*{\em Proof of the monotonicity property~(\ref{mon_salpha})}
For~\(\alpha \leq 1\) we couple the original Poisson process~\({\cal P}\) and the homogenous process~\({\cal P}_{\delta}\) in the following way. Let~\(V_{i}, i \geq 1\) be i.i.d.\ random variables each with density~\(f(.)\) and let~\(N_V\) be a Poisson random variable with mean~\(n,\) independent of~\(\{V_i\}.\) The nodes~\(\{V_i\}_{1 \leq i \leq N_V}\) form a Poisson process with intensity~\(nf(.)\) which we denote as~\({\cal P}\) and colour green.

Let~\(U_i, i \geq 1\) be i.i.d.\ random variables each with density~\(\epsilon_2-f(.)\) where~\(\epsilon_2 \geq 1\) is as in~(\ref{f_eq}) and let~\(N_U\) be a Poisson random variable with mean~\(n(\epsilon_2-1).\) The random variables~\((\{U_i\},N_U)\) are independent of~\((\{V_i\},N_V)\) and the nodes~\(\{U_i\}_{1 \leq i  \leq N_U}\) form a Poisson process with intensity~\(n(\epsilon_2-f(.))\) which we denote as~\({\cal P}_{ext}\) and colour red. The nodes of~\({\cal P}\) and~\({\cal P}_{ext}\) together form a homogenous Poisson process with intensity~\(n\epsilon_2,\) which we denote as~\({\cal P}_{\delta}\) and define it on the probability space~\((\Omega_{\delta},{\cal F}_{\delta}, \mathbb{P}_{\delta}).\)

Let~\(\omega_{\delta} \in \Omega_{\delta}\) be any configuration and let~\({\cal I} := \{i_{j}\}_{1 \leq j \leq Q}\) be the indices of the dense squares in~\(\{R_j\},\) each containing at least three nodes of~\({\cal P}\) and let~\({\cal I}_{\delta} := \{i^{(\delta)}_{j}\}_{1 \leq j \leq Q_{\delta}}\) be the indices of the~\(\delta-\)dense squares in~\(\{R_j\},\) each containing at least three nodes of~\({\cal P}_{\delta}.\) We set~\({\cal I} = \emptyset\) if~\(\omega_{\delta}\) contains no dense square and set~\({\cal I}_{\delta} = \emptyset\) if~\(\omega_{\delta}\) contains no~\(\delta-\)dense square. By definition~\({\cal I} \subseteq {\cal I}_{\delta} \) and defining~\(S_{\alpha} = S_{\alpha}(\omega_{\delta})\) and~\(S^{(\delta)}_{\alpha} = S_{\alpha}^{(\delta)}(\omega_{\delta})\) as before, we have that~\(S_{\alpha}\) is determined only by the green nodes of~\(\omega_{\delta}\) while~\(S^{(\delta)}_{\alpha}\) is determined by both green and red nodes of~\(\omega_{\delta}.\)

From the monotonicity property, we therefore have that~\(S_{\alpha}(\omega_{\delta}) \leq S^{(\delta)}_{\alpha}(\omega_{\delta})\) and so for any~\(x > 0\) we have
\begin{equation}\label{mon_eq}
\mathbb{P}_{\delta}(S^{(\delta)}_{\alpha} < x) \leq \mathbb{P}_{\delta}(S_{\alpha} < x)  = \mathbb{P}_0(S_{\alpha} < x),
\end{equation}
proving~(\ref{mon_salpha}).

If~\(\alpha > 1,\) we perform a slightly different analysis. Letting~\(\epsilon_1 \leq 1\) be as in~(\ref{f_eq}), we construct a Poisson process~\({\cal P}_{ext}\) with intensity~\(n(f(.)-\epsilon_1)\) and colour nodes of~\({\cal P}_{ext}\) red. Letting~\({\cal P}_{\delta}\) be another independent Poisson process with intensity~\(n\epsilon_1,\) we colour nodes of~\({\cal P}_{\delta}\) green. The superposition of~\({\cal P}_{ext}\) and~\({\cal P}_{\delta}\) is a Poisson process with intensity~\(nf(.),\) which we define on the probability space~\((\Omega_{\delta},{\cal F}_{\delta},\mathbb{P}_{\delta}).\) In this case, the sum~\(S_{\alpha}\) is determined by both green and red nodes while~\(S^{(\delta)}_{\alpha}\) is determined only by the green nodes. Again using the monotonicity property of~\(S_{\alpha},\) we get~(\ref{mon_eq}).~\(\qed\)

\subsection*{\em Additive relations}
If~\(TSP(x_1,\ldots,x_j), j \geq 1\) denotes the length of the TSP cycle with vertex set~\(\{x_1,\ldots,x_j\},\)
then for any~\(k \geq 1,\) we have
\begin{equation}\label{mono_2}
TSP(x_1,\ldots,x_{j+k}) \leq TSP(x_1,\ldots,x_j) + TSP(x_{j+1},\ldots,x_{j+k}) + (c_2\sqrt{2})^{\alpha}
\end{equation}
and if~\( \alpha \leq 1\) and the edge weight function~\(h\) is a metric, then
\begin{equation}
TSP(x_1,\ldots,x_j) \leq TSP(x_1,\ldots,x_{j+k}). \label{mono}
\end{equation}
\emph{Proof of~(\ref{mono_2}) and~(\ref{mono})}: To prove~(\ref{mono_2}),
suppose~\({\cal C}_{1}\) is the minimum weight spanning cycle formed
by the nodes~\(\{x_l\}_{1 \leq l \leq j}\) and~\({\cal C}_2\) is the minimum weight spanning cycle
formed by~\(\{x_l\}_{j+1 \leq l \leq j+k}.\) Let~\(e_1 = (u_1,v_1) \in {\cal C}_1\)
and~\(e_2 = (u_2,v_2) \in {\cal C}_2\) be any two edges. The cycle
\[{\cal C}_{tot} = \left({\cal C}_1 \setminus \{e_1\}\right) \cup \left({\cal C}_2 \setminus \{e_2\}\right) \cup \{(u_1,u_2), (v_1,v_2)\}\]
obtained by removing the edges~\(e_1,e_2\) and adding the ``cross" edges~\((u_1,u_2)\) and~\((v_1,v_2)\) is a spanning cycle containing all the nodes~\(\{x_l\}_{1 \leq l \leq j+k}.\) The edges~\((u_1,u_2)\) and~\((v_1,v_2)\) have a Euclidean length of at most~\(\sqrt{2}\)
and so a weight of at most~\((c_2\sqrt{2})^{\alpha}\) using the bounds for the metric~\(h\) in~(\ref{eq_met1}).
This proves~(\ref{mono_2}).

It suffices to prove~(\ref{mono}) for~\(k =1.\) Let~\({\cal C} = (y_1,\ldots,y_{j+1},y_1)\) be any cycle
with vertex set~\(\{y_i\}_{1 \leq i \leq j+1} = \{x_i\}_{1 \leq i \leq j+1}\)
and without loss of generality suppose that~\(y_{j+1} = x_{j+1}.\)
Removing the edges~\((y_j,y_{j+1})\)
and~\((y_{j+1},y_1),\) and adding the edge~\((y_1,y_j)\)
we get a new cycle~\({\cal C}'\) with vertex set~\(\{x_i\}_{1 \leq i \leq j}.\)

Since the edge weight function~\(h\) is a metric, we have by triangle inequality that~\(h(y_j,y_1) \leq h(y_j,y_{j+1}) + h(y_{j+1},y_1).\)
Using~\((a+b)^{\alpha} \leq a^{\alpha} + b^{\alpha}\) for~\(a,b > 0\) and~\(0 < \alpha \leq 1\) we get
that~\(h^{\alpha}(y_j,y_1) \leq h^{\alpha}(y_j,y_{j+1}) + h^{\alpha}(y_{j+1},y_1).\)
Therefore the weight~\(W({\cal C}')\) of~\({\cal C}'\)
\[W({\cal C'}) = \sum_{i=1}^{j-1}h^{\alpha}(y_i,y_{i+1}) + h^{\alpha}(y_j,y_1)
\leq \sum_{i=1}^{j} h^{\alpha}(y_i,y_{i+1}) + h^{\alpha}(y_{j+1},y_1) = W({\cal C}).\]
Therefore~\(TSP(x_1,\ldots,x_j) \leq W({\cal C}') \leq W({\cal C}).\)
Taking minimum over all cycles~\({\cal C}\)
with vertex set~\(\{x_i\}_{1 \leq i \leq j+1},\) we get~(\ref{mono}) for~\(k=1.\)~\(\qed\)

\subsection*{\em Moments of random variables}
Let~\(X \geq 1\) be any integer valued random variable such that
\begin{equation}\label{x_dist}
\mathbb{P}(X \geq l) \leq e^{-\theta (l-1)}
\end{equation}
for all integers~\(l \geq 1\) and some constant~\(\theta > 0\) not depending on~\(l.\) For every integer~\(r \geq 1,\)
\begin{equation}\label{disc_tel}
\mathbb{E}X^{r} \leq r\sum_{l\geq 1}  l^{r-1} \mathbb{P}(X \geq l) \leq r\sum_{l \geq 1} l^{r-1} e^{-\theta (l-1)} \leq \frac{r!}{(1-e^{-\theta})^{r}}
\end{equation}
\emph{Proof of~(\ref{disc_tel})}: For~\(r \geq 1\) we have
\begin{equation}\label{tel_sum}
\mathbb{E}X^{r} = \sum_{l \geq 1} l^{r} \mathbb{P}(X= l) = \sum_{l \geq 1}l^{r} \mathbb{P}(X \geq l) - l^{r}\mathbb{P}(X \geq l+1)
\end{equation}
and substituting the~\(l^{r}\) in the final term of~(\ref{tel_sum}) with~\(  (l+1)^{r} - ((l+1)^{r}-l^{r})\) we get
\begin{eqnarray}
\mathbb{E}X^r &=& \sum_{l \geq 1} \left(l^{r} \mathbb{P}(X \geq l) - (l+1)^{r} \mathbb{P}(X \geq l+1)\right)  \nonumber\\
&&\;\;\;\;\;\;\;\; + \;\;\;\sum_{l \geq 1} ((l+1)^{r}-l^{r}) \mathbb{P}(X \geq l+1) \nonumber\\
&=& 1 + \sum_{l \geq 1}((l+1)^{r}-l^{r}) \mathbb{P}(X \geq l+1)  \nonumber\\
&=& \sum_{l \geq 0}((l+1)^{r}-l^{r}) \mathbb{P}(X \geq l+1)  \label{gent}
\end{eqnarray}
where the second equality is true since~\(l^{r}\mathbb{P}(X \geq l) \leq l^{r}e^{-\theta (l-1)} \longrightarrow 0\) as~\(l~\rightarrow~\infty.\) Using~\((l+1)^{r} - l^{r} \leq r\cdot (l+1)^{r-1}\) in~(\ref{gent}), we get the first relation in~(\ref{disc_tel}).

We prove the second relation in~(\ref{disc_tel}) by induction as follows. Let~\(\gamma = e^{-\theta} < 1\) and~\(J_r := \sum_{l \geq 1} l^{r-1} \gamma^{l-1}\) so that
\begin{equation}
J_{r+1}(1 - \gamma) = \sum_{l \geq 1} l^{r} \gamma^{l-1} - \sum_{l \geq 1}l^{r} \gamma^{l} = \sum_{l \geq 1} \left(l^{r}-(l-1)^{r} \right)\gamma^{l-1}. \nonumber
\end{equation}
Using~\(l^{r}-(l-1)^r \leq r\cdot l^{r-1}\) for~\(l \geq 1\) we therefore get that
\[J_{r+1}(1-\gamma) \leq r\sum_{l \geq 1}l^{r-1} \gamma^{l-1} = rJ_r\]
and so the second relation in~(\ref{disc_tel}) follows from induction.~\(\qed\)

\subsection*{\em Scaling and translation relations}
For a set of nodes~\(\{x_1,\ldots,x_n\}\) in the unit square~\(S,\)
recall from Section~\ref{intro} that~\(K_n(x_1,\ldots,x_n)\)
is the complete graph formed by joining all the nodes by straight line segments and the edge~\((x_i,x_j)\)
is assigned a weight of~\(d^{\alpha}(x_i,x_j),\)
where~\(d(x_i,x_j)\) is the Euclidean length of the edge~\((x_i,x_j).\)
We denote~\(TSP(x_1,\ldots,x_n)\) to be the length of the minimum spanning cycle of~\(K_n(x_1,\ldots,x_n)\)
with edge weights obtained as in~(\ref{min_weight_cycle}),

\emph{\underline{Scaling}}: For any~\(a > 0,\) consider the graph~\(K_n(ax_1,\ldots,ax_n)\) where the length of the edge between
the vertices~\(ax_1\) and~\(ax_2\) is simply~\(a\) times the length of the edge between~\(x_1\) and~\(x_2\)
in the graph~\(K_n(x_1,\ldots,x_n)\)
Using the definition of TSP in~(\ref{min_weight_cycle}) we then have~\(TSP(ax_1,\ldots,ax_n) = a^{\alpha} TSP(x_1,\ldots,x_n)\)
and so if~\(Y_1,\ldots,Y_n\) are~\(n\) nodes uniformly distributed in the square~\(aS\) of side length~\(a,\)
then
\begin{equation}\nonumber
TSP(n;a) := TSP(Y_1,\ldots,Y_n) = a^{\alpha} TSP(X_1,\ldots,X_n),
\end{equation}
where~\(X_i = \frac{Y_i}{a}, 1 \leq i \leq n\) are i.i.d.\ uniformly distributed in~\(S.\)
Recalling the notation~\(TSP_n = TSP(X_1,\ldots,X_n)\) from~(\ref{min_weight_cycle}) we therefore get
\begin{equation}\label{tsp_scale}
\mathbb{E} TSP(n;a) = a^{\alpha} \mathbb{E}TSP_n.
\end{equation}

\emph{\underline{Translation}}: For~\(b \in \mathbb{R}^2\) consider the graph~\(K_n(x_1+b,\ldots,x_n+b).\) Using the translation property~\((b2),\)
the weight~\(h(x_1+b,x_2+b) \leq h_0 \cdot h(x_1,x_2),\) the weight of the edge between~\(x_1\) and~\(x_2.\)
Using the definition of TSP in~(\ref{min_weight_cycle})
we therefore have
\begin{equation}\label{trans_rel}
TSP(x_1+b,\ldots,x_n+b) \leq h_0^{\alpha} \cdot TSP(x_1,\ldots,x_n),
\end{equation}
obtaining the desired bound.~\(\qed\)

\subsection*{Acknowledgement}
I thank Professors Rahul Roy, Thomas Mountford, Federico Camia, C. R. Subramanian and the referee for crucial comments that led to an improvement of the paper. I also thank IMSc for my fellowships.

\bibliographystyle{plain}

\begin{thebibliography}{10}





\bibitem{beard} Beardwood, J., Halton, J. H. and Hammersley, J. M. (1959)
\newblock{The shortest path through many points}.
\newblock{\em Proceedings Cambridge Philosophical Society}, \textbf{55}, pp. 299--327.

\bibitem{brad} Bradonjic, M., Elsasser, R., Friedrich, T., Sauerwald, T. and Stauffer, A.  (2010)
\newblock{Efficient Broadcast on Random Geometric Graphs}.
\newblock{\em Proceedings SODA 2010}, pp. 1412--1421.


\bibitem{gutin} Gutin, G. and Punnen, A. P. (2006)
\newblock{\em The traveling salesman problem and its variations}.
\newblock{Springer}.

\bibitem{steele2} Steele, J. M. (1981)
\newblock{Subadditive Euclidean functionals and nonlinear growth in geometric probability}.
\newblock{\em Annals of Probability}, \textbf{9}, pp. 365--376.

\bibitem{steele3} Steele, J. M. (1997)
\newblock{\em Probability Theory and Combinatorial Optimization}.
\newblock{SIAM}.

\bibitem{stin} Steinerberger, S. (2015)
\newblock{New bounds for the Traveling Salesman constant}.
\newblock{\em Advances in Applied Probability}, \textbf{47}, pp. 27--36.

\bibitem{rhee} Rhee, W. T.  (1991)
\newblock{On the Fluctuations of the Stochastic Traveling Salesperson Problem}.
\newblock{\em Mathematics of Operation Research}, \textbf{16}, pp. 482--489.

\bibitem{yuk} Yukich, J. (1998)
\newblock{\em Probability Theory of Classical Euclidean Optimization Problems}.
\newblock{Lecture Notes in Mathematics}, \textbf{1675}, Springer.





\end{thebibliography}

\end{document}